\documentclass{article}

     \PassOptionsToPackage{numbers}{natbib}
     \usepackage[preprint]{neurips_2022}

\usepackage[utf8]{inputenc} % allow utf-8 input
\usepackage[T1]{fontenc}    % use 8-bit T1 fonts
\usepackage{hyperref}       % hyperlinks
\usepackage{url}            % simple URL typesetting
\usepackage{booktabs}       % professional-quality tables
\usepackage{amsfonts}       % blackboard math symbols
\usepackage{nicefrac}       % compact symbols for 1/2, etc.
\usepackage{microtype}      % microtypography
\usepackage{xcolor}   

\usepackage{amsmath}
\usepackage{amssymb}
\usepackage{graphicx}
\usepackage{multirow}
\usepackage{enumerate}
\usepackage{arydshln}
\usepackage{multicol}
\usepackage{algorithm}
\usepackage{algpseudocode}
\usepackage[utf8]{inputenc} % allow utf-8 input
\usepackage[T1]{fontenc}    % use 8-bit T1 fonts
\usepackage{hyperref}       % hyperlinks

\usepackage{amsthm}

\usepackage{url}            % simple URL typesetting
\usepackage{booktabs}       % professional-quality tables
\usepackage{amsfonts}       % blackboard math symbols
\usepackage{nicefrac}       % compact symbols for 1/2, etc.
\usepackage{microtype}      % microtypography
\usepackage{xcolor}         % colors
\usepackage{wrapfig}
\usepackage{subcaption}
\usepackage{multicol}
\def\x{\mathbf{x}}
\def\y{\mathbf{y}}
\def\z{\mathbf{z}}

\def\v{\mathbf{v}}

\def\I{\mathbf{I}}

\def\S{\mathcal{S}}

\def\A{\mathbf{A}}

\newtheorem{theorem}{Theorem}[section]
\newtheorem{proposition}{Proposition}[section]

\newtheorem{exam}{Example}
\newtheorem{assumption}{Assumption}[section]

\title{Towards Extremely Fast Bilevel Optimization with Self-governed Convergence Guarantees}

\author{%
	Risheng~Liu$^{\ast 1}\quad $ Xuan~Liu$^{1}\quad$   Wei~Yao$^{2}\quad$  Shangzhi~Zeng$^{3}\quad$ Jin~Zhang$^{2,4}$\\
	$^{1}$International School of Information Science \& Engineering, DUT\\
	$^{2}$Department of Mathematics, SUSTech\quad$^{3}$Department of Mathematics and Statistics, UVIC \\
	$^{4}$National Center for Applied Mathematics Shenzhen\\
}

\begin{document}

\maketitle

\begin{abstract}
Gradient methods have become mainstream techniques for Bi-Level Optimization (BLO) in learning and vision fields. The validity of existing works heavily relies on solving a series of approximation subproblems with extraordinarily high accuracy. Unfortunately, to achieve the approximation accuracy requires executing a large quantity of time-consuming iterations and computational burden is naturally caused. This paper is thus devoted to address this critical computational issue. 
In particular, we propose a single-level formulation to uniformly understand existing explicit and implicit Gradient-based BLOs (GBLOs). This together with our designed counter-example can clearly illustrate the fundamental numerical and theoretical issues of GBLOs and their naive accelerations. By introducing the dual multipliers as a new variable, we then establish Bilevel Alternating Gradient with Dual Correction (BAGDC), a general framework, which significantly accelerates different categories of existing methods by taking specific settings. A striking feature of our convergence result is that, compared to those original unaccelerated GBLO versions, the fast BAGDC admits a unified non-asymptotic convergence theory towards stationarity. A variety of numerical experiments have also been conducted to demonstrate the superiority of the proposed algorithmic framework.
\end{abstract}

\section{Introduction}
Bi-Level Optimization (BLO) has recently attracted growing interests due to its wide applications in the field of deep learning \cite{liu2021investigating}, especially for hyperparameter optimization \cite{franceschi2017forward,okuno2018hyperparameter,mackay2018self}, meta learning~\cite{franceschi2018bilevel,rajeswaran2019meta,zugner2018adversarial}, neural architecture search~\cite{liu2018darts,liang2019darts+,chen2019progressive}, adversarial learning~\cite{pfau2016connecting}, and reinforcement learning~\cite{yang2019provably}, etc. BLO tackles nested optimization structures appearing in applications. In the last decade, BLO has emerged as a prevailing optimization technique for modern machine learning tasks with an underlying hierarchy. 

Recently, Gradient-based BLOs (GBLOs) have become mainstream techniques for BLO in learning and vision fields. In this type, an approximate version of BLO can be solved by taking into explicit account the optimization dynamics for the Lower-Level (LL) objective. The issue of the quality of this approximation has been well studied in \cite{franceschi2017forward,liu2020generic,grazzi2020iteration,ji2021bilevel}.  A major advantage of GBLOs is that an approximate hypergradient can be directly calculated by automatic differentiation based on the trajectory of the LL variable. However, GBLOs theoretically require a large number of LL optimization steps to achieve the desired accuracy \cite{liu2020generic}. Unfortunately, to achieve the accuracy requires executing a large quantity of time-consuming iterations. Thus, computational and memory burden are naturally caused. In practice, one simple heuristic to accelerate is to only perform a small number of iterations for the LL problem, see, e.g., neural architecture search~\cite{liu2018darts,elsken2020meta}, adversarial learning~\cite{goodfellow2014generative} and reinforcement learning~\cite{pfau2016connecting}. However, the theoretical properties are still unclear, although competitive empirical performance has been reported. 

%Recently, a broad collection of algorithms have been proposed to solve (stochastic) BLOs with using the one-step LL problem iteration idea, see, e.g., \cite{chen2022single,khanduri2021near,yang2021provably}. However, a large quantity of time-consuming matrix multiplications or Hessian-vector products is required there. Moreover, the strong convexity of the LL objective is required for obtaining the convergence results of all these proposed algorithms.

On the other hand, to remove the need for differentiating through the LL optimization path, especially when the LL dynamic system iterates many times, a type of implicit gradient-based BLO methods (called I-GBLOs for simplicity) is employed in~\cite{pedregosa2016hyperparameter,rajeswaran2019meta,lorraine2020optimizing} for hyperparameter optimization and meta-learning. In fact, I-GBLOs first leverages the implicit differentiation to derive an analytical expression of hypergradient, and then solves inexactly the LL problem up to a tolerance and performs an approximate matrix inversion accurately. Hence, I-GBLOs requires executing a large quantity of iterations and computational burden is naturally caused.

\subsection{Our Motivations and Contributions}\label{Our Motivations and Contributions}
As mentioned above, to access the hypergradient accurately, the number of iteration steps in both GBLOs and I-GBLOs theoretically goes to infinity \cite{liu2020generic,pedregosa2016hyperparameter,rajeswaran2019meta}. Hence, computational burden is naturally caused. One simple heuristic to accelerate is to only performs one-step iteration for the LL problem, see, e.g., DARTS \cite{liu2018darts}  and its applications for various tasks, cf. e.g., neural architecture search~\cite{liu2018darts,elsken2020meta}, adversarial learning~\cite{goodfellow2014generative}, reinforcement learning~\cite{pfau2016connecting}, hyperparameter optimization~ \cite{cao2021pancsc} , meta learning~\cite{finn2017model} and image processing~\cite{jin2021bridging}. Although competitive empirical performance has been reported widely, the corresponding theoretical properties are still unclear. Recently, a broad collection of algorithms have been proposed to solve (stochastic) BLOs with using the one-step LL problem iteration idea, see, e.g., \cite{chen2022single,khanduri2021near,yang2021provably}. However, a large quantity of time-consuming matrix multiplications or Hessian-vector products is required there. Moreover, the strong convexity of the LL objective is required for obtaining the convergence results of all these proposed algorithms.

In this paper, we first construct an interesting example in Section~\ref{3.1} explicitly to indicate that existing approaches including GBLOs and I-GBLOs with adopted naive acceleration may lead to incorrect solutions. Therefore, an important question which motives this work arises.
{\it Can we design a extremely fast bilevel  acceleration strategy which enjoy favorably both extremely low computational cost and reliable theoretical guarantee?}
%{\it Can we design a bilevel acceleration strategy which enjoy favorably both extremely low computational cost and theoretical convergence guarantee?}

%{\it Can we design a extremely fast bilevel  acceleration strategy which enjoy favorably both extremely low computational cost and theoretical convergence guarantee (including convergence rate and convergence theory towards statinary) ?}
In response to the above question, we first replace the LL problem with its first-order optimality condition and the resulting problem is a single-level optimization problem with equality constraints.. Thus we can introduce naturally the dual multipliers as a new variable. Our first insight is to well exploit the bilevel structure, keeping the update of the LL variable by gradient descent of the LL objective and its variants. Our second insight is that a prospective point where the hypergradient of BLO vanishes happens to be a KKT solution of the above single-level optimization problem. Thereby, to guarantee the non-asymptotic convergence, we update cautiously the Upper-Level (UL) and dual multipliers variables with the LL variable,  alternatively according to the KKT condition. Compared with naive accelerations, the newly involved multipliers variable can be regarded as a dual correction. Therefore, we establish Bilevel Alternating Gradient with Dual Correction (BAGDC), an extremely fast BLO algorithm, to ease both computational and memory burden. Furthermore, compared to those original unaccelerated GBLO and I-GBLO versions, we analyze the non-asymptotic convergence behavior of BAGDC towards self-governed KKT stationary solutions.

%In response to the above question, we first replace the LL problem with its first-order optimality condition and the resulting problem is a single-level optimization problem with equality constraints. Thus we can introduce naturally the dual multipliers as a new variable. Our first insight is that, the equality constraints come from the LL minimization problem of BLO. Hence, to well exploit the bilevel structure, it is reasonable for us to keep the update of the LL variable by gradient descent of the LL objective and its variants.  Our second insight is that a prospective point where the hypergradient of BLO vanishes happens to be a KKT solution of the above single-level optimization problem. 
%Hence, to guarantee the non-asymptotic convergence, we update cautiously the Upper-Level (UL) and dual multipliers variables with the LL variable by alternating gradient, according to the form of the KKT condition. Compared with naive accelerations, the newly involved multipliers variable can be regarded as a dual correction. Therefore, we establish Bilevel Alternating Gradient with Dual Correction (BAGDC), an extremely fast BLO algorithm, to ease both computational and memory burden. Furthermore, compared to those original unaccelerated GBLO and I-GBLO versions, we analyze the non-asymptotic convergence behavior of BAGDC towards KKT stationary solutions. 

The main contributions of this paper can be summarized as follows:
\begin{itemize}
	%\item We first design an interesting counter-example to illustrate such invalidation of the naive accelerations of the existing Gradient-based BLOs (GBLOs) and I-GBLOs. Then we propose a single-level formulation of BLO to uniformly understand existing explicit and implicit GBLOs. 
	
	\item We first demonstrate fundamental theoretical issues arising from the naive single-inner-step acceleration which has been widely-witnessed in many real-world learning and vision applications, by explicitly exhibiting the trajectory of GBLOs on our designed counter-example. 
	
	%\item We first investigate the trajectory of GBLOs on our designed counter-example to demonstrate that there actually exist fundamental theoretical issues in these naive single-inner-step acceleration for mainstream GBLOs, although which have been widely applied in many real-world learning and vision applications.
	
	%\item By introducing the dual multipliers as a new variable, we establish Bilevel Alternating Gradient with Dual Correction (BAGDC), an extremely fast BLO algorithm, in which the LL, the dual multipliers and UL variables are updated by alternating gradient. Hence, BAGDC admits a much simpler implementation. As comparison, the existing BLO methods including EG and ig always update the LL and UL variables in a nested manner.
	
	\item By reformulating BLO as a single-level optimization with equality constraints and introducing dual multipliers from the viewpoint of KKT condition, we then establish BAGDC, an alternating gradient method, to successfully avoid solving a series of time-consuming nested optimization subproblems, and hence remarkably faster than a variety of existing GBLOs and I-GBLOs.
	
	%\item By reformulating BLO as a single-level optimization with equality constraints and introducing dual multipliers, from the viewpoint of KKT condition, we then establish BAGDC, an alternating gradient method, that is extremely faster than a variety of existing GBLOs and I-GBLOs (all require to solve a series of time-consuming nested optimization subproblems).
	
	%\item By reformulating BLO as a single-level optimization with equality constraints and  introducing the dual multipliers, from the view point of the KKT condition, we then establish BAGDC, an alternating gradient method, that is extremely faster than a variety of existing GBLOs and I-GBLOs (all require to solve a series of time-consuming nested optimization subproblems).
	
%	\item Compared to those original unaccelerated GBLO versions, the fast BAGDC admits a unified non-asymptotic convergence theory towards KKT stationary points. It is worth pointing out that our result is also applicable to BLOs with multiple LL minimal points. 
	
%	\item More importantly, rather than these relatively weak convergence results (i.e., without convergence rate or convergence theory towards stationary) in exiting works, we prove a series of new theories to guarantee the non-asymptotic convergence of BAGDC towards KKT stationary points. It is worth pointing out that our result is also applicable to BLOs with multiple LL minimal points.

\item An important feature of our theoretical result relies on a series of new non-asymptotic convergence guarantees of BAGDC towards KKT stationary points.
Compared to exiting works which are either lack of convergence towards stationarity or convergence rate,
our result is evidently stronger and general enough to cover BLOs with multiple LL minimizers.
\end{itemize}

%We first investigate the trajectory of GBLOs on our designed counter-example to demonstrate that although have been widely applied in many real-world learning and vision applications, there actually exist fundamental theoretical issues in these naïve single-inner-step acceleration for mainstream GBLOs.
%
%By introducing a Lagrange-based formulation to analyze the optimality conditions of BLOs, we then establish BAGDC, an alternating gradient method that is extremely faster than a variety of existing explicit and implicit GBLOs, which all require a series of time-consuming nested optimization subproblems.
%
%More importantly, rather than these relatively weak convergence results (i.e., without convergence rate and KKT guarantees) in exiting works, we prove a series of new theories to guarantee the non-asymptotic convergence of BAGDC towards KKT stationary points. It is worth pointing out that our result is also applicable to BLOs with multiple LL minimal points.

\section{Bilevel Alternating Gradient with Dual Correction} \label{sec2}

In this work, we consider the BLO problem in the following form:
\begin{equation}\label{blo problem}
	\mathop{\min}\limits_{\x, \y} F(\x,\y), \quad \mathrm{s.t.} \quad  \y \in \S(\x):=\arg\min_{\y \in \mathbb{R}^m} f(\x,\y),
\end{equation}
where $\x\in\mathbb{R}^n$, $\y\in\mathbb{R}^m$ are UL and LL variables respectively. Here the UL and LL objectives are smooth functions with Lipschitz continuous first and second order derivatives.\footnote{This assumption is common in the BLO literature, see, e.g., \cite{ghadimi2018approximation,ji2021bilevel,khanduri2021near,chen2022single}. A formal description of the assumption is stated in the Supplemental Material.} 

\subsection{Fundamental Issues and Counter-Example} \label{3.1}

BLOs are usually considered as nested structural optimization problems. Thus, the existing BLO methods including GBLOs and I-GBLOs always update the UL and LL variables in a nested manner. This roughly causes the computational and theoretical issues for the existing methods. To design an extremely fast BLO algorithm, we aim to update the UL and LL variables by alternating gradient, cf. the well-known Gradient Descent-Ascent (GDA) for minimax optimization problems. But alternating gradient for BLOs is challenging since naive accelerations of the existing approaches may lead to incorrect solutions. Let us take a widely used acceleration technique to exemplify this issue.  Indeed, if one only performs one-step LL iteration in GBLOs, the resulting algorithm, called Naive One-Step Acceleration (NOSA for simplicity of notation), reads as 
\begin{equation}\label{darts}
		\begin{aligned}
			\y_{k+1}=&\y_k-\beta\nabla_{\y} f(\x_{k},\y_k),\\
			\x_{k+1}=&\x_k-\alpha\left(\nabla_{\x} F(\x_k,\y_{k+1})-\beta\nabla_{\x\y}^2 f(\x_k,\y_{k})\nabla_{\y} F(\x_k,\y_{k+1})\right),
	\end{aligned}
\end{equation}
where $\alpha,\beta$ are stepsizes. Next, we design an interesting counter-example (Example~\ref{counter} below) to illustrates that NOSA in Eq.~(\ref{darts}) may lead to incorrect solutions. 
\begin{exam}\label{counter} (Counter-Example) Consider the following BLO problem:
	\begin{eqnarray}\label{eq:counter}
		\min _{\x \in \mathbb{R}^n} \frac{1}{2}\|\x-\z_0\|^2+\frac{1}{2}\y^*(\x)^\top\A \y^*(\x) \quad \mathrm { s.t. } \quad\y^*(\x) = \arg \min _{\y \in \mathbb{R}^n} f(\x, \y)=  \frac{1}{2}\y^\top\A\y-\x^\top \y,
	\end{eqnarray} 
	where $\x\in \mathbb{R}^n$, $\y\in \mathbb{R}^n$, $\A\in\mathbb{S}^{n\times n}$ is a positive-definite symmetric matrix and $\z_0\neq0$ is a given point in $\mathbb{R}^n$. Note that the LL objective is strongly convex w.r.t. the LL variable and $\x,\y$ have the same dimension. Assuming that $\left(\x_{k}, \y_{k}\right)$ in Eq.~\eqref{darts} converges to $(\bar{\x}, \bar{\y})$. Then $\A\bar{\y}-\bar{\x}=0$ and $(\bar{\x}-\z_0)+\beta \A\bar{\y}=0$. This implies that $\bar{\x}=\z_0/(1+\beta)$. Now we show that $\bar{\x}$ is not a stationary point of $\varphi(\x):=F(\x,\y^*(\x))$ if $\z_0$ is not a eigenvector of $\A$. Indeed, a simple computation gives $\nabla \varphi(\bar{\x})=(\bar{\x}-\z_0)+\A^{-1} \bar{\x}=\left(\A^{-1} \z_0-\beta \z_0\right)/(1+\beta)\neq0$. 
\end{exam}
By using the same counter-example~\ref{counter} and a similar argument, one can easily check that the same result also holds for other naive one-step accelerations of I-GBLOs, such as Naive One-Step Accelerations of Conjugate Gradient (CG) and Neumann series (NS) methods. For more details, we refer the reader to the Supplemental Material.

\subsection{A Single-Level Reformulation and Dual Correction}

To figure out the reason why NOSA in Eq.~(\ref{darts}) does not work well, we first transform BLO into a single-level optimization problem with equality constraint:
\begin{equation}\label{ecp}
	\min_{\x,\y}\  F(\x,\y)\quad
	\mathrm{s.t.}\quad \nabla_{\y} f(\x,\y)=0.
\end{equation} 
Clearly, BLO in Problem~\eqref{blo problem} and Problem~\eqref{ecp} are equivalent when the LL objective $f(\x,\y)$ is convex w.r.t. $\y$ for any fixed $\x$. Denote $\mathcal{L}(\x,\y,\v)=F(\x,\y)-\v^T\nabla_{\y} f(\x,\y)$, where $\v\in\mathbb{R}^m$ is the dual multipliers. We say that $(\x^*,\y^*)$ is a KKT point of Problem~\eqref{ecp} if there exists a dual multipliers $\v^*\in\mathbb{R}^m$ such that the following KKT condition holds: 
\begin{align}\label{KKT}
	\nabla\mathcal{L}(\x^*,\y^*,\v^*)
	=\left(
	\begin{array}{ccc}
		\nabla_{\x} F(\x^*,\y^*)-\nabla_{\x\y}^2 f(\x^*,\y^*)\v^*\\
		\nabla_{\y} F(\x^*,\y^*)-\nabla_{\y\y}^2 f(\x^*,\y^*)\v^*\\
		-\nabla_{\y} f(\x^*,\y^*)
	\end{array}
	\right)=0.
\end{align}
To better describe what we are going to do with the KKT condition, we briefly provide a strong connection between the KKT condition and the hypergradient using in the existing BLO methods. In fact, the first part (called $\x$-counterpart) of the KKT condition is in consistence with the hypergradient, which is equal to $\nabla_{\x} F(\x^*,\y^*)-\nabla_{\x\y}^2 f(\x^*,\y^*)\v^*$ when both the second part (called $\y$-counterpart) and the third one (called feasible condition) of the KKT condition hold, if the Hessian $\nabla_{\y\y}^2 f(\x^*,\y^*)$ is nonsingular. The latter condition is common in the literature. But the nonsingularity of $\nabla_{\y\y}^2 f(\x^*,\y^*)$ is not necessary for the KKT condition. Therefore, the KKT condition is more applicable to general BLOs, e.g., BLOs with multiple LL minimal points. 

%It is worth noting that $\v^*=\left[\nabla_{\y\y}^2 f(\x^*,\y^*)\right]^{-1} \nabla_{\y} F(\x^*,\y^*)$ 

Now we can uniformly understand the existing explicit and implicit Gradient-based BLOs. From the viewpoint of non-asymptotic convergence analysis, roughly speaking, the existing BLO methods can be divided into three steps: (1) given UL variable, update the LL variable by solving the LL optimization problem up to tolerance, which satisfies prescribed dynamic accuracy requirements. From the view point of the KKT condition, this step is to guarantee the feasible condition; (2) given UL and LL variable, update the dual multipliers by approximating the correct dual multipliers accurately. There are mainly two kinds of techniques, i.e., automatic differentiation based on the trajectory of LL variable in RHG~\cite{franceschi2017forward} and solving a linear system based on the Conjugate Gradient (CG) or Neumann Series (NS)~\cite{rajeswaran2019meta,lorraine2020optimizing,grazzi2020iteration}. The correct dual multipliers here is the one satisfying the $\y$-counterpart of the KKT condition; (3) update the UL variable by the accurate approximate hypergradient. As we discussed in the previous paragraph, this corresponds to the $\x$-counterpart of the KKT condition. From the above analysis, one can clearly see the reason why NOSA in Eq.~\eqref{darts} does not work well. In fact, NOSA did not update the dual multipliers properly. The hidden update $\v_{k+1}= \beta\nabla_{\y} F(\x_k,\y_{k+1})$ in NOSA can not guarantee the $\y$-counterpart of the KKT condition.

The above issue motivate us to consider the following question: how to update the dual multipliers if the update of the LL variable by gradient descent of the LL objective is reserved? This question is important and nontrivial since the dual multipliers are always updated by the equality constraint in the standard methods of solving nonlinear optimization with equality constraints, see, e.g.,~\cite[Chapter17]{wright1999numerical}. The reason why we prefer updating the LL variable by $\y_{k+1}=\y_k-\beta\nabla_{\y} f(\x_{k},\y_k)$ and its variants is that the equality constraint of Problem~\eqref{ecp} comes from the LL minimization problem of BLO. Hence, to well exploit the bilevel structure, it is reasonable for us to keep the update of the LL variable by gradient descent. This can also guarantee the feasible condition at the limit point. We also keep the update of the UL variable by $\x_{k+1}=\x_k-\alpha\left(\nabla_{\x} F(\x_k,\y_{k+1})-\left[\nabla_{\x\y}^2 f(\x_k,\y_{k})\right]\v_{k+1}\right)$ and its variants to ensure the $\x$-counterpart of the KKT condition. Hence, to guarantee the convergence of algorithm towards self-governed KKT stationary solutions, the update of the multipliers variable should be to ensure the $\y$-counterpart of the KKT condition. 

When the LL problem is merely convex, BLO can be with multiple LL minimizers. Then updating of the LL variable by gradient descent of the LL objective may not necessarily generate a sequence towards minimizer of UL objective over LL problem solutions. Inspired by BDA \cite{liu2020generic}, we use aggregation parameters to aggregate both UL and LL objectives. Specifically, we define 
	\begin{equation}
		\psi_{\mu}(\x,\y):=\mu \lambda F(\x,\y)+(1-\mu)f(\x,\y),
	\end{equation}
	where $\mu,\lambda$ are aggregation parameters. We then replace the LL objective $f$ by $\psi_{\mu_k}$ and let the aggregation parameter $\mu_k$ goes to $0$ as the iteration goes on. When the LL problem is strongly convex w.r.t. $\y$, we just take the aggregation parameter $\mu=0$ in the following algorithm.
Now, we propose the following Bilevel Alternating Gradient with Dual Correction (BAGDC) in Algorithm~\ref{alg:BAGDC} for BLOs.

\begin{algorithm}[h]
	\caption{Bilevel Alternating Gradient with Dual Correction (BAGDC)}\label{alg:BAGDC}
	\hspace*{0.02in} {\bf Initialize:} %算法的输入， \hspace*{0.02in}用来控制位置，同时利用 \\ 进行换行
	$\x_0,\y_0,\v_0$, stepsizes $\alpha_k,\beta_k,\eta_k$, aggregation parameters $\mu_k,\lambda$;\\
	\hspace*{0.02in} {\bf Let:}  $\psi_{\mu_k}(\x,\y):=\mu_k \lambda F(\x,\y)+(1-\mu_k)f(\x,\y)$;
	%	\hspace*{0.02in} {\bf Output:} %算法的结果输出
	%	output result
	%	\hspace*{0.02in} {\bf Set:} aggregation function $\psi_{\mu,\lambda}(\x,\y)=\mu\lambda F(\x,\y)+(1-\mu) f(\x,\y)$;
	\begin{algorithmic}[1]
		%		\State some description % \State 后写一般语句
		\For{$k=0,1,\dots,K-1$} % For 语句，需要和EndFor对应
		\State update $\y_{k+1}=\y_k-\beta_k\nabla_{\y}\psi_{\mu_k}(\x_k,\y_k)$;
		\State update $\v_{k+1}=\v_k+\eta_k\left(\nabla_{\y} F(\x_k,\y_{k+1})- \left[\nabla_{\y\y}^2 \psi_{\mu_k}(\x_k,\y_{k+1})\right]\v_k\right);$
		\State update 
		$\x_{k+1}=\x_k-\alpha_k\left(\nabla_{\x} F(\x_k,\y_{k+1})-\left[\nabla_{\x\y}^2 \psi_{\mu_k}(\x_k,\y_{k})\right]\v_{k+1}\right).$
		%		\If{condition} % If 语句，需要和EndIf对应
		%		\State ...
		%		\Else
		%		\State ...
		%		\EndIf
		\EndFor
		%		\While{condition} % While语句，需要和EndWhile对应
		%		\State ...
		%		\EndWhile
		%		\State \Return result
	\end{algorithmic}
\end{algorithm}

To help understand BAGDC in Algorithm~\ref{alg:BAGDC}, the following points should be emphasized.
\begin{itemize}
	\item In BAGDC, to handle a general BLO, e.g., BLOs with multiple LL minimizers, we use aggregation parameters $\mu_k,\lambda$ to aggregate the LL objective $f$ and a positive multiple of the UL objective $\lambda F$. An explicitly strategy for choosing stepsizes $\alpha_k,\beta_k,\eta_k$ and aggregation parameters $\mu_k,\lambda$ is presented in Section~\ref{Convergence Analysis} based on theoretical convergence analysis. Specially, when the LL objective is strongly convex w.r.t. $\y$, the aggregation parameter $\mu$ is taken to be $0$. In this case, the stepsizes $\alpha,\beta,\eta$ can be taken to be constant stepsizes satisfying certain appropriate conditions given in Section~\ref{Convergence Analysis}. 
	
	\item Compared with naive accelerations, e.g., NOSA in Eq.~\eqref{darts}, the newly involved multipliers variable $\v$ can be regarded as a dual correction. Its appearance turns bilevel alternating gradient into a practical acceleration strategy which enjoy favorably both extremely simpler implementation and theoretical convergence guarantee. For more details we refer the reader to the theoretical investigation in Section~\ref{Convergence Analysis} and experiments in Section~\ref{experiments}.
\end{itemize}

\section{Theoretical Investigation}\label{Convergence Analysis}

This section is devoted to the convergence analysis of our proposed algorithm BAGDC with and without strong convexity assumption of LL problem. The UL and LL objectives $F(\x,\y)$ and $f(\x,\y)$ are assumed to be smooth functions with Lipschitz continuous first and second order derivatives.\footnote{This assumption is common in the BLO literature. A formal description of the assumption is stated in the Supplemental Material.}  Please notice that all the proofs of our theoretical results given here and other auxiliary results are stated in the Supplemental Material.

\subsection{General Results for LL Merely Convex Case}

%\begin{assumption}\label{assump}
%	$F(\x,\y)$ and $f(\x,\y)$ are smooth functions with Lipschitz continuous first and second order derivatives.\footnote{ A formal description of the assumption is stated in the Supplemental Material.} 
%\end{assumption}

To handle general BLOs, e.g., BLOs with multiple LL minimizers, we take the following standard assumption as in \cite{liu2020generic}.
\begin{assumption}\label{assump_convex}
	For any fixed $\x$, $f(\x,\cdot)$ is convex and $F(\x,\cdot)$ is $\sigma_F$-strongly convex.
\end{assumption}
Recall that 
	\begin{equation*}
	\psi_{\mu}(\x,\y):=\mu \lambda F(\x,\y)+(1-\mu)f(\x,\y),
\end{equation*}
where $\mu,\lambda$ are aggregation parameters.
Let $0<\mu\leq 1/2$ and $\lambda>0$, for any fixed $\x$, as $F(\x,\cdot)$ is assumed to be $\sigma_F$-strongly convex, the aggregation function $\psi_{\mu}(\x,\y)$ is $\sigma_{\psi_{\mu}}$-strongly convex in $\y$ with $\sigma_{\psi_{\mu}}=\mu\lambda\sigma_{F} > 0$. Hence $\psi_{\mu}(\x,\y)$ has a unique minimal point for every $\x$, denoted by $\y^*_{\mu}(\x)$.
Under the smoothness assumption of $f$ and $F$, $\psi_{\mu}(\x,\y)$ is also $L_{\psi_{\mu}}$-smooth in $\y$ for any fixed $\x$.

To analyze the convergence of BAGDC, we define the following Lyapunov function for the generated sequence $\{(\x_k,\y_k,\v_k)\}$ by BAGDC~\ref{alg:BAGDC}:
\begin{equation}
	V_k:=F(\x_k,\y^*_{\mu_k}(\x_k))+\frac{1}{2}\|\y_k-\y^*_{\mu_k}(\x_k)\|^2+\frac{1}{2}\|\v_k-\v^*_{\mu_k}(\x_k)\|^2,
\end{equation}
where $\v^*_{\mu}(\x)=\left[\nabla_{\y\y}^2 \psi_{\mu}(\x,\y^*_{\mu}(\x))\right]^{-1}\nabla_{\y} F(\x,\y^*_{\mu}(\x))$ is the correct dual multipliers. The first term quantifies the descent of the overall UL objective functions, the second term characterizes the LL solution errors, and the third term characterizes the multipliers errors. We can prove the following descent property of Lyapunov function $V_k$ for BAGDC.

\begin{proposition}\label{propVk}
	Suppose Assumption \ref{assump_convex} holds. If we choose aggregation parameter $\mu_k \downarrow 0$ and step sizes
	\begin{equation}\label{stepsize1}
		\beta_k\leq\frac{2}{\sigma_{\psi_{\mu_{k}}}+L_{\psi_{\mu_{k}}}}, \quad \eta_k \leq  c_\eta \mu_k^4 \beta_k,\quad \alpha_{k}\leq c_\alpha \mu_{k}^7 \eta_{k},
	\end{equation}
	then when $\mu_k$ is sufficiently small, the sequence of $\x_k,\y_k,\v_k$ generated by BAGDC satisfies
	\begin{align}\label{diffV}
		V_{k+1}-V_k\leq & -\frac{\alpha_k}{2}\|\nabla \varphi_{\mu_{k}}(\x_k)\|^2-\frac{\widehat{\alpha}_k}{2}\|\x_{k+1}-\x_k\|^2	
		-\frac{\widehat{\beta}_k}{2}\|\y_k-\y_{\mu_k}^*(\x_k)\|^2-\frac{\widehat{\eta}_k}{2}\|\v_k-\v_{\mu_k}^*(\x_k)\|^2\nonumber\\
		& + C_1\left(\frac{\mu_k-\mu_{k+1}}{\mu_{k}}\right) + C_2 \frac{1}{\alpha_{k}}\left(\mu_{k}-\mu_{k+1}\right)^2\nonumber + C_3\frac{1}{\alpha_k}\left(\frac{\mu_{k}}{\mu_{k+1}}\right)^2\left(\mu_{k}-\mu_{k+1}\right)^2,
	\end{align}
	where $\widehat{\alpha}_k$, $\widehat{\beta}_k$, $\widehat{\eta}_k$, $c_\eta$, $c_\alpha$ $C_1$, $C_2$ and $C_3$ are all positive constants, whose formal formulae are stated in the  Supplemental Material.
\end{proposition}

According to this descent property, we can notice that to guarantee that $\|\nabla \varphi_{\mu_{k}}(\x_k)\|$  converges to $0$, we need to choose aggregation parameter $\mu_k$ and step size $\alpha_k$ appropriately such that $$\sum_{k=0}^\infty \alpha_k=\infty, \ ~\sum_{k=0}^\infty \frac{1}{\alpha_k}\left(\mu_k-\mu_{k+1}\right)^2<\infty,\ \ \mathrm{and}\ ~\sum_{k=0}^\infty \frac{1}{\alpha_k}\left(\frac{\mu_{k}}{\mu_{k+1}}\right)^2\left(\mu_k-\mu_{k+1}\right)^2<\infty.$$

In following, we present a simple and handy strategy for choosing aggregation parameter $\mu_k$ and step sizes $\alpha_k$, $\beta_k$ and $\eta_k$ for guaranteeing the convergence of our proposed BAGDC.

\begin{theorem}\label{thmconvex}
	Suppose Assumption \ref{assump_convex} holds. If we choose $\mu_k= \bar{\mu} \left(\frac{1}{k+1}\right)^p$ with $0<p<\frac{1}{11}$ and 
	\begin{equation}
		\alpha_k= \bar{\alpha} \mu_{k}^{11},\quad \beta_k=\bar{\beta}, \quad \eta_k= \bar{\eta} \mu_{k}^{4},
	\end{equation}
	where $\bar{\beta} \in (0, 1/L_{\psi_{\mu_{k}}})$,
	then when all of positive constants $\bar{\alpha}$, $\bar{\eta}$, $\bar{\mu}$ are sufficiently small, the sequence of $\x_k,\y_k$ generated by BAGDC satisfies
%	\begin{equation*}
%		\min_{1\leq k\leq K}\left\{\|\nabla \varphi_{\mu_{k}}(\x_k)\|^2+\frac{\widehat{\beta}_k}{\alpha_{k}}\|\y_k-\y^*_{\mu_k}(\x_k)\|^2
%		+\frac{\widehat{\eta}_k}{\alpha_{k}}\|\v_k-\v^*_{\mu_k}(\x_k)\|^2\right\}=O\left(\frac{\ln K}{K^{1-11p}}\right).
%	\end{equation*}
%	\begin{equation*}
%	\min_{1\leq k\leq K}\left\{\|\nabla \varphi_{\mu_{k}}(\x_k)\|^2\right\}=O\left(\frac{\ln K}{K^{1-11p}}\right)\quad\mathrm{and}\quad
%	\min_{1\leq k\leq K}\left\{\|\y_k-\y^*_{\mu_k}(\x_k)\|^2
%	\right\}=O\left(\frac{\ln K}{K^{1-11p}}\right).
%\end{equation*}
\begin{equation*}
	\min_{0\leq k\leq K}\left\{\|\nabla \varphi_{\mu_{k}}(\x_k)\|^2\right\}=O\left(\frac{\ln K}{K^{1-11p}}\right)\quad\mathrm{and}\quad
	\min_{0\leq k\leq K}\left\{\|\y_k-\y^*_{\mu_k}(\x_k)\|^2
	\right\}=O\left(\frac{\ln K}{K^{1-p}}\right).
\end{equation*}
%	\begin{equation*}
%	\min_{1 \le k \le K} \{\|\nabla \varphi_{\mu_{k}}(\x_k)\| \}= O\left(\sqrt{\frac{\ln K}{K^{1-11p}}}\right),\ \mathrm{and}\ \min_{1 \le k \le K} \{ \|\y_k-\y^*_{\mu_k}(\x_k)\| \} = O\left(\sqrt{\frac{\ln K}{K}}\right).	
%	\end{equation*}
\end{theorem}

%To appropriately characterize the convergence rate of BAGDC, we define the KKT residual for $\x$ as following
%\[
%\mathrm{KKT}(\x, \mu):= \|\nabla\mathcal{L}(\x,\y^*_{\mu}(\x), \v^*_{\mu}(\x))\|.
%\]
%Then, $\mathrm{KKT}(\x, \mu) = 0$ if and only if $(\x,\y^*_{\mu}(\x))$ is a KKT point to Eq.~\eqref{ecp} and $\v^*_{\mu}(\x)$ is the multipliers.
%Under the linear independent constraint qualification (LICQ), we can have following non-asymptotic convergence rate result of BAGDC from the perspective of KKT residual.

To appropriately characterize the convergence rate of BAGDC, we define the KKT residual for $(\x,\y,\v)$ as following
\[
\mathrm{KKT}(\x, \y,\v):= \|\nabla \mathcal{L}(\x,\y, \v)\|^2.
\]
Then, $\mathrm{KKT}(\x, \y,\v) = 0$ if and only if $(\x,\y)$ is a KKT point to Problem~\eqref{ecp} and $\v$ is a multipliers.
Under the linear independent constraint qualification (LICQ), we can have following non-asymptotic convergence rate result of BAGDC from the perspective of KKT residual.

\begin{theorem}
	Suppose Assumption \ref{assump_convex} holds and aggregation parameter $\mu_k$ and step sizes $\alpha_k$, $\beta_k$ and $\eta_k$ are chosen as same as in Theorem \ref{thmconvex}.
	Let $(\x_k,\y_k,\v_k)$ be the sequence of generated by BAGDC with sufficiently small $\bar{\alpha}$, $\bar{\eta}$, $\bar{\mu}$, if $(\x_k,\y_k)$ is bounded and for any limit point $(\bar{\x}, \bar{\y})$ of $(\x_k,\y_k)$, LICQ of Problem~\eqref{ecp} holds at $(\bar{\x}, \bar{\y})$, i.e., rows of $\left(\nabla_{\x\y}^2f(\bar{\x}, \bar{\y})^T, \nabla_{\y\y}^2f(\bar{\x}, \bar{\y})^T\right)$ are linearly independent, then we have
	\[
	\min_{0 \leq k \leq K} \left\{ \mathrm{KKT}(\x_k, \y_k,\v_k) \right\}= O\left(\max\left\{ \frac{\ln K}{K^{1-11p}}, \frac{1}{K^{2p}}\right\}\right).
	\]
%	\[
%	\min_{1 \le k \le K} \left\{ \mathrm{KKT}(\x_k, \mu_k) \right\}= O\left(\max\left\{ \sqrt{\frac{\ln K}{K^{1-11p}}}, \left(\frac{1}{K+1} \right)^p \right\}\right).
%	\]
	Furthermore, any limit point $(\bar{\x}, \bar{\y})$ of $\left(\x_k,\y_k\right)$ is a KKT point of Problem~\eqref{ecp}.
\end{theorem}

\subsection{Special Discussion for LL Strongly Convex Case}

For the special case where the LL objective $f(\x,\cdot)$ is strongly convex with constant $\sigma_{f}>0$ for any fixed $\x$, the strong convexity assumption on UL objective $F$ is not needed, and BAGDC can be shown to have a $O(1/K)$ non-asymptotic convergence rate with taking $\mu_k=0$ and constant step sizes $\alpha_k$, $\beta_k$ and $\eta_k$.
We use the following assumption to replace Assumption \ref{assump_convex} in this section.
\begin{assumption}\label{assump_convex2}
	For any fixed $\x$, $f(\x,\cdot)$ is $\sigma_{f}$-strongly convex.
\end{assumption}

Note that when the LL objective $f(\x,\cdot)$ is strongly convex, $\y^*(\x)$ is well defined and be differentiable. Hence, it is proper to use the gradient of $\varphi(\x)=F(\x,\y^*(\x))$ to characterize the optimality of the generated sequence by the algorithms. 
\begin{theorem}\label{theorem_Strong_convex}
	Suppose Assumption \ref{assump_convex2} holds. If we choose $\mu_k=0$ and 
	\begin{equation}
		\alpha_k= \bar{\alpha},\quad \beta_k=\bar{\beta}, \quad \eta_k = \bar{\eta},
	\end{equation}
	where $\bar{\alpha}, \bar{\beta}$, and $\bar{\eta}$ are positive constants satisfying certain conditions \footnote{A formal description of the conditions is stated in the Supplemental Material.},
	then $\x_k$ generated by BAGDC satisfies
	\[
	\min_{0 \le k \le K} \{\|\nabla \varphi(\x_k)\|^2 \}= O\left(\frac{1}{K}\right).
	\]
\end{theorem}

\section{Experiments}\label{experiments}

In this section, we conduct extensive experiments on both synthetic data and real applications to verify our theoretical findings and demonstrate the efficacy of BAGDC, especially on very high-dimensional tasks.
%In this section, we conduct extensive experiments on both synthetic data and real applications to very our theoretical findings and demonstrate the efficacy of BAGDC, especially on very high-dimensional tasks. 
The experiments were done on a PC with Intel i7-9700K CPU (4.2 GHz), 32GB RAM, and an NVIDIA RTX 2060S GPU. The algorithm was implemented using PyTorch 1.6 with CUDA 10.2.

\subsection{Numerical Verification}\label{4.1}

In this subsection, we start by discussing the effect of dual correction under different settings. We then compare the performance of BAGDC with other GBLOs and I-GBLOs by using different accuracies for LL problem. The third section demonstrates the efficacy of BAGDC compared to RHG, CG and NS. Finally, we compare the empirical performance of BAGDC for BLOs with multiple LL minimizers.
We conduct experiments in Eq.~\eqref{eq:counter} which satisfies the assumptions in Theorem~\ref{theorem_Strong_convex}. 
In particular we set $\A=\I$, $\z_0=\mathbf{e}$, where $\mathbf{e}$ is the vector whose elements are all equal to 1.  The unique solution is $\x^{*}=\y^{*}=\mathbf{e}/2$. For the existing method, we set the number of lower level iterations to 100.  $\mathbf{d}=\mathbf{d}(\x,\y,\v)$ is the direction of descent of $\x$  calculated by different methods, e.g. for NOSA, $\mathbf{d} = \left(\nabla_{\x} F(\x,\y)-\left[\nabla_{\x\y}^2 f(\x,\y)\right]\nabla_{\y}F(\x,\y)\right)$. We define the gradient of $\x$  as $\nabla \varphi$ where $\varphi(\x)=F(\x,\y^*(\x))$. 
%Note that we use $\Vert\mathbf{d}\Vert$ to indicate that the algorithm converges to its stationary point, while the rest indicator containing the distance to the optimum point to indicate whether the algorithm converges to the optimum point.

\textbf{Effect of Dual Correction with Different Settings.}
Although we proved in Section~\ref{Convergence Analysis} that the addition of dual variable $\v$ guarantees convergence, the update step size $\eta$ of $\v$ still depends on properties of $F$ and $f$ (e.g. the Lipschitz constant and the strongly convex coefficients). Since this is usually agnostic in practical applications, we first verify the corrective effect of $\v$ for convergence for different $\eta$ settings. In addition to different fixed $\eta$, we also validate a strategy for the Adaptive (Adapt.) computation of $\eta$, a CG-inspired strategy.\footnote{More details is stated in the Supplemental Material}.
Figure~\ref{fig:eta} (a) shows the convergence of NOSA and BAGDC with different $\eta$. 
It can be seen that all $\eta$ improve the convergence performance of non-convergence NOSA. For a fixed $\eta$, any $\eta$ within the theoretical bound is guaranteed to converge, although smaller $\eta$  lead to a lower speed of convergence. On the other hand, $\eta$ outside the theoretical bound may lead to a collapse after a while. Note that despite the convergence result of Theorem \ref{theorem_Strong_convex}, BAGDC is not guaranteed to be monotonically decreasing, and in fact, some degree of oscillation is expected for some $\eta$. $\eta$ with adaptive strategy achieved a convergence speed second only to $\eta=C_f$. Since adaptive $\eta$ does not need to know the nature of $F$ and $f$, we used adaptive $\eta$ as our strategy in all subsequent experiments.
However, considering that adaptive method may require more computational cost in computing $\eta$, we further show in Figure~\ref{fig:eta} (b) the Average (Avg.) step time and the overall time to converge to $\Vert\mathbf{d}\Vert\leq 10^{-5}$ for different $\eta$ settings. As can be seen, although the adaptive method takes slightly more time for a single step, the overall convergence time is still second only to $\eta=C_f$. 
%
%
%we can observe the convergence of $\v$ is guaranteed for $\eta$ within the theoretical interval, where the theoretical upper bound converges fastest. Smaller $\eta$ also guarantees convergence but leads to lower speeds and possible oscillations. Adaptive $\eta$, on the other hand, maintains a sub-optimal convergence rate in the absence of oscillations. Figure~\ref{fig:eta} (b) shows  the average time per step and the total convergence time for the different methods to reach convergence.  Although the adaptive $\eta$ leads to slightly higher computation times each step, it is still second only to $C_f$ for the total convergence time. Since the $C_f$ may not be computable in real applications, 

\begin{figure}[!htbp]
	\centering 
	\setlength{\tabcolsep}{1.4mm}{
		\begin{tabular}{c@{\extracolsep{0.35em}}c@{\extracolsep{0.35em}}c@{\extracolsep{0.35em}}c@{\extracolsep{0.35em}}c@{\extracolsep{0.35em}}c@{\extracolsep{0.35em}}}
			\multicolumn{4}{c}{\includegraphics[height=0.014\textheight]{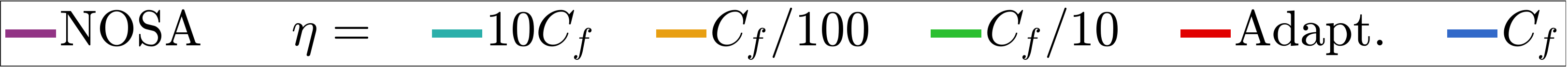}}&& \\
			{\includegraphics[width=0.18\linewidth,trim=0 -18 0 0,clip]{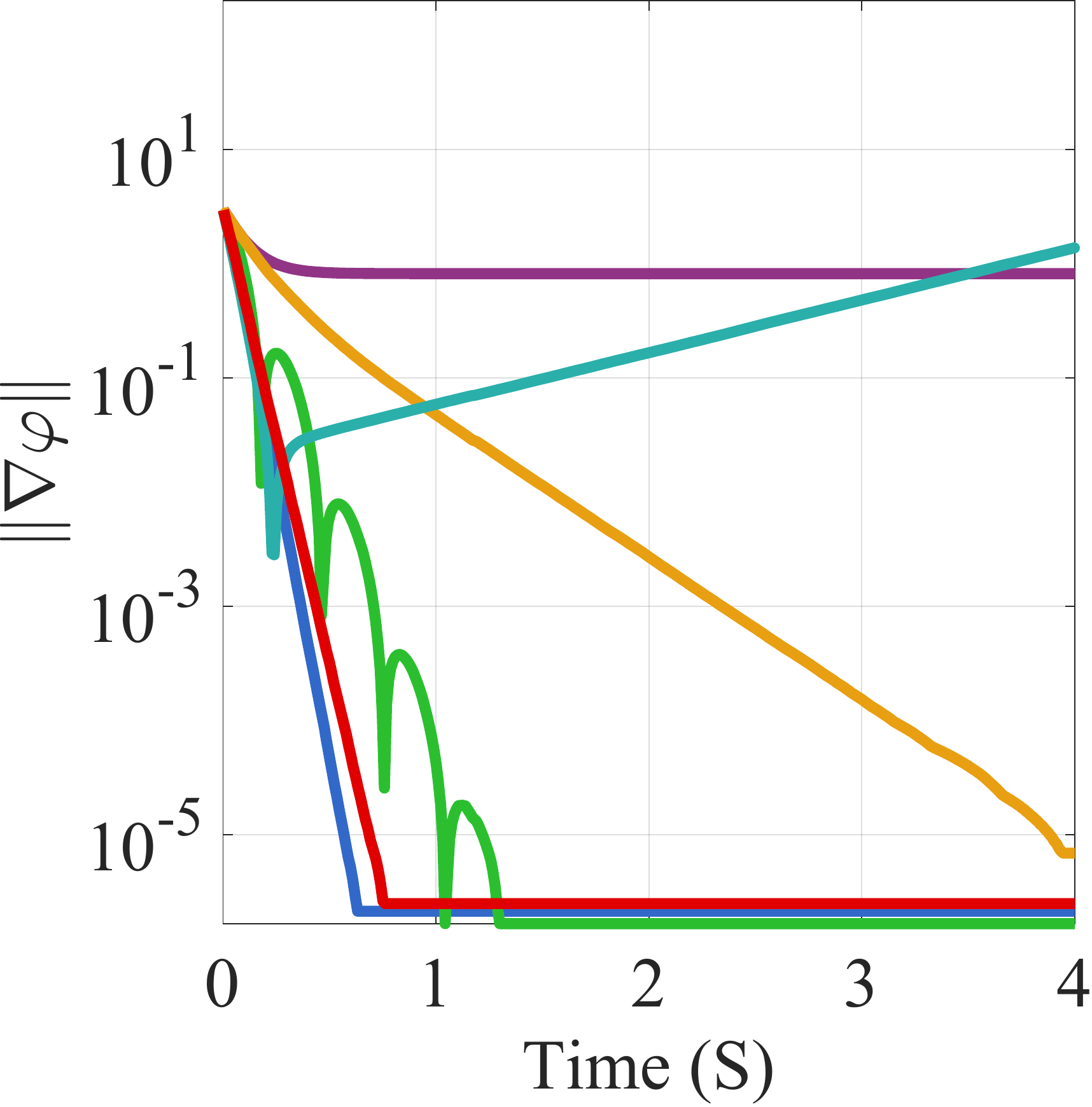}}&
			{\includegraphics[width=0.18\linewidth,trim=0 -18 0 0,clip]{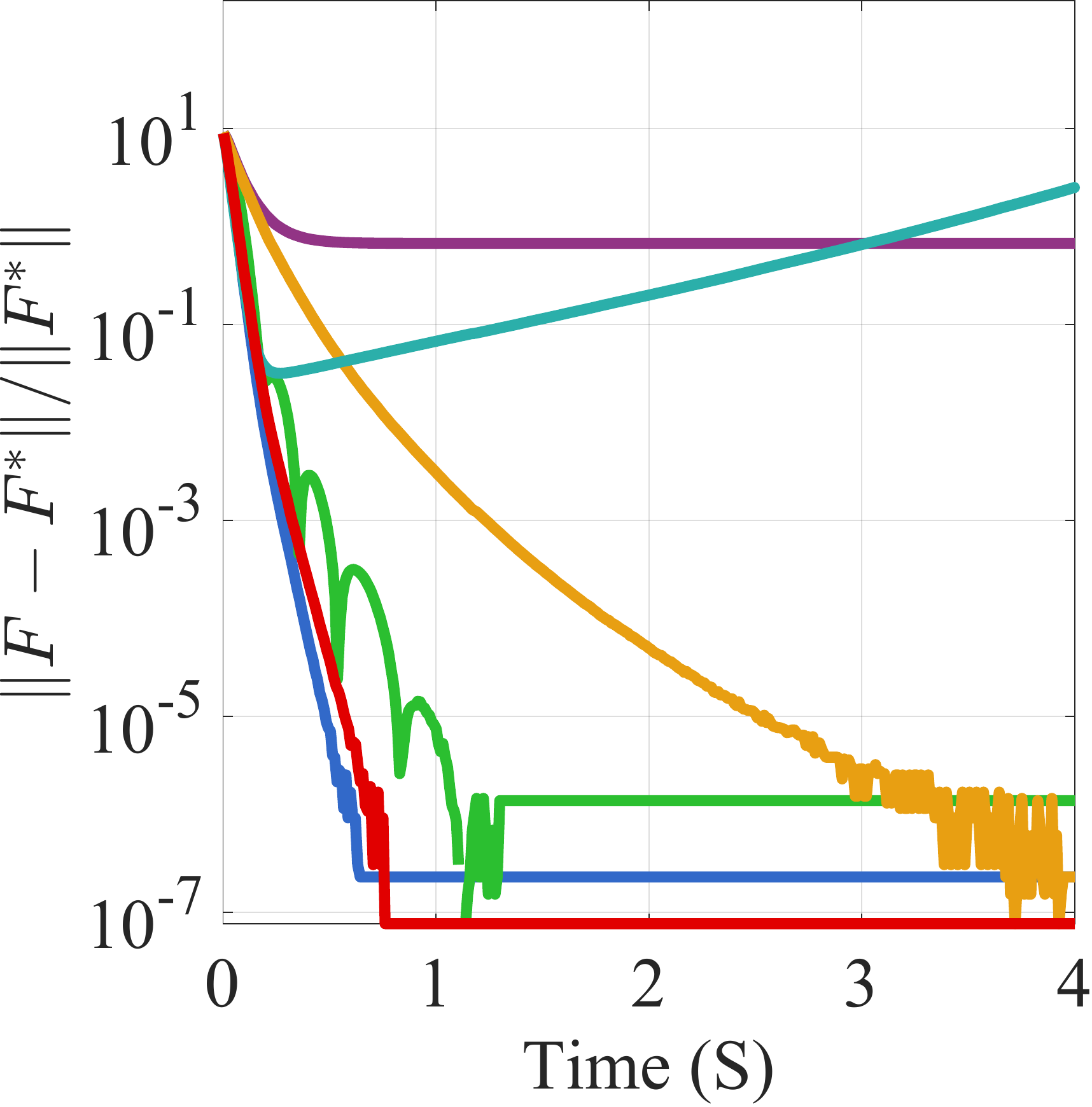}}&
			{\includegraphics[width=0.18\linewidth,trim=0 -18 0 0,clip]{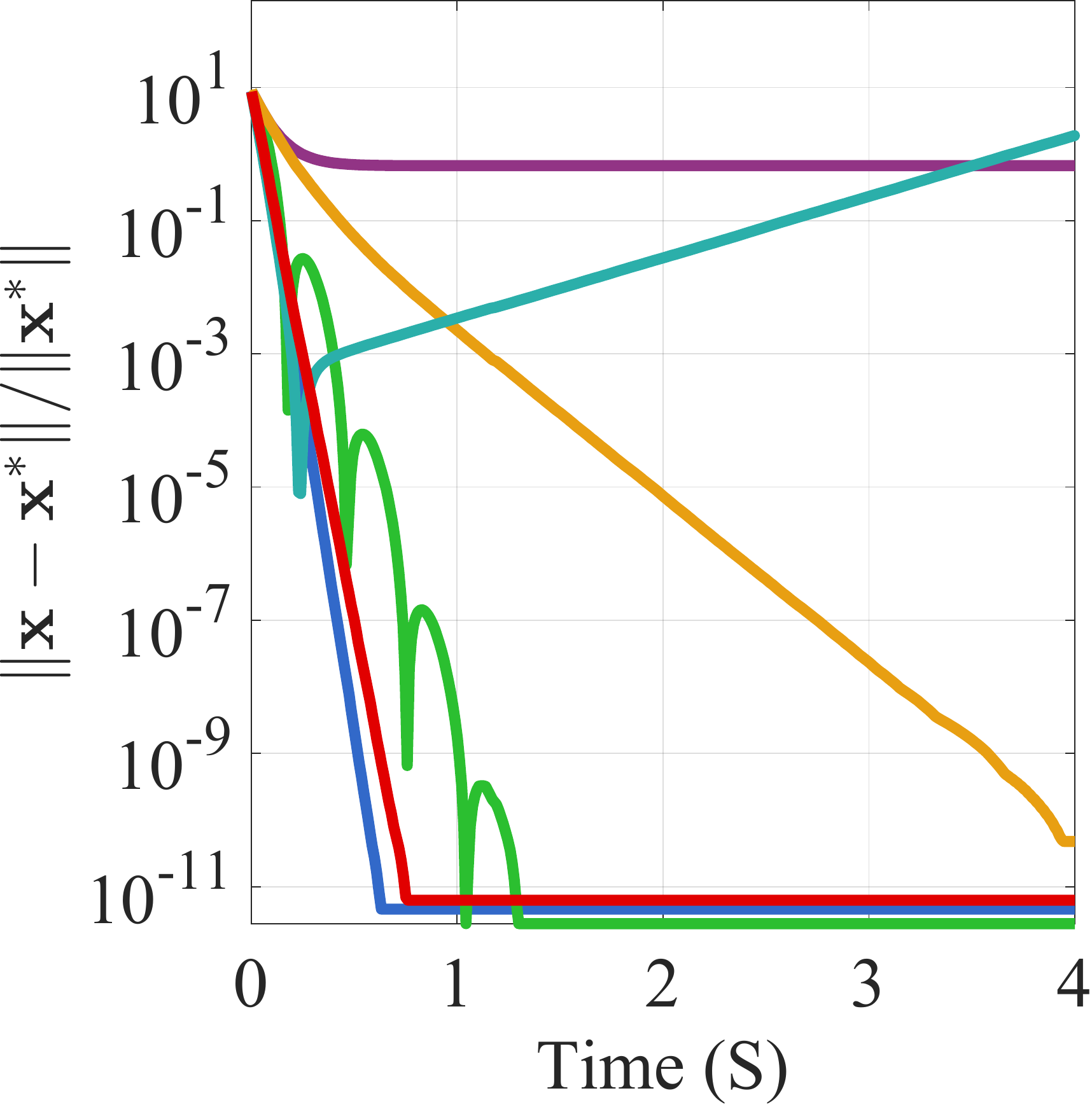}}&
			{\includegraphics[width=0.18\linewidth,trim=0 -18 0 0,clip]{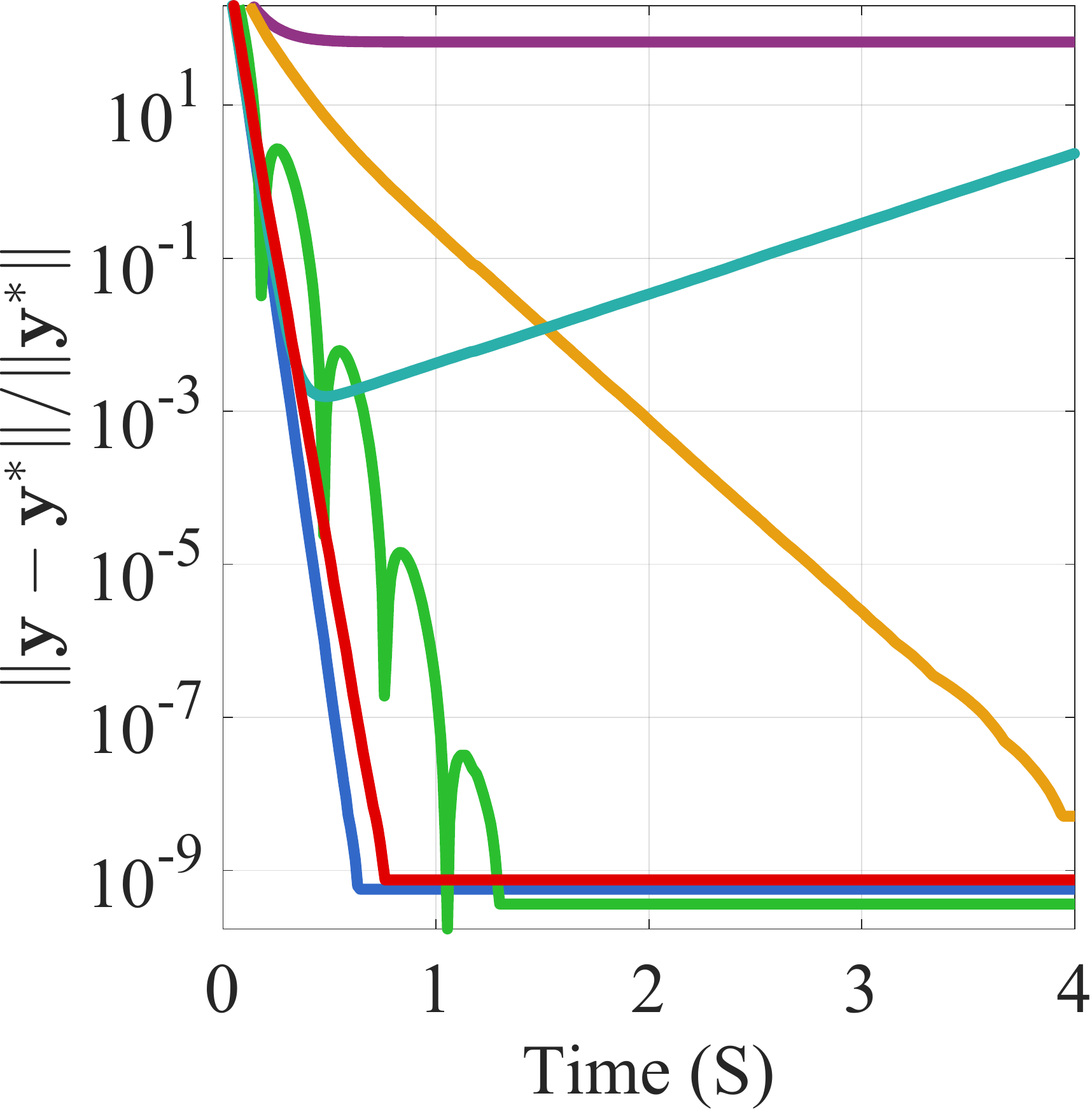}}&&
			{\includegraphics[width=0.20\linewidth,trim=0 3 0 0,clip]{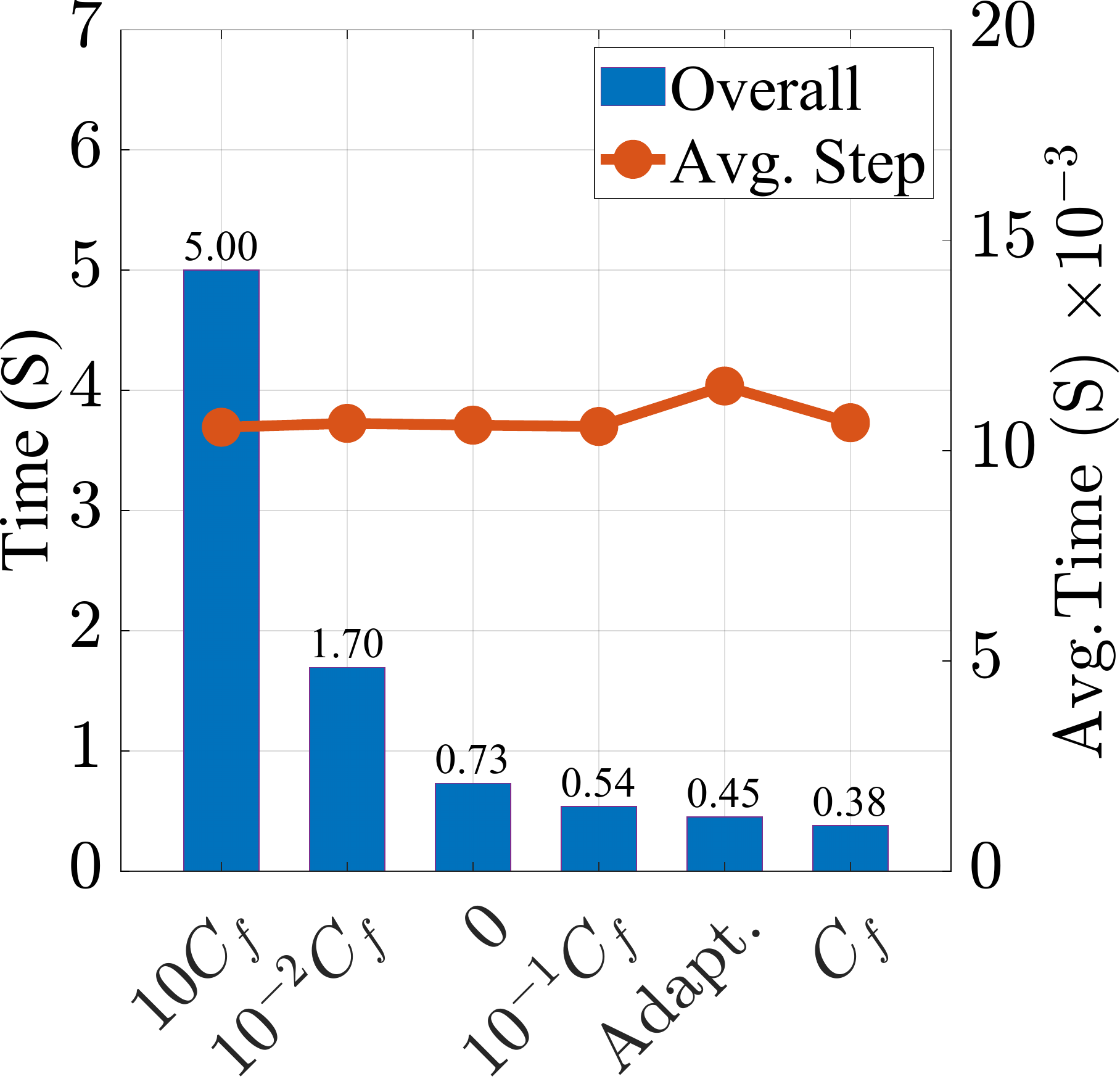}}\\
			%\multicolumn{6}{c}{\includegraphics[height=0.015\textheight]{fig/toyeta/legend.pdf}} \\
			\multicolumn{4}{c}{\footnotesize(a) Convergence curves in different $\eta$ }& \multicolumn{2}{c}{\footnotesize(b) Overall/Avg. Step}\\
	\end{tabular}}
	\caption{Convergence comparison of update strategies for $\v$ in different $\eta$. $C_f$ is the theoretical upper bound. (a) Convergence curve with different  $\eta$. (b) Overall time for convergence and Average (Avg.) step time.
		The theoretical upper bound has the best convergence speed, but the sub-optimal adaptive $\eta$ can be applied to real applications where the theoretical upper bound is difficult to compute.
	}\label{fig:eta}
	\vspace{-0.2cm}
\end{figure}

%Since BAGDC use alternating optimization, we only need to solve the LL problem once at each iteration. Other methods that use nested optimization require the exact solution of the LL problem at each iteration.

\begin{figure}[!htbp]
	%	\vspace{-0.1cm} 
	\setlength{\tabcolsep}{1.4mm}{
		\begin{tabular}{l@{\extracolsep{0.03em}}l@{\extracolsep{0.03em}}l@{\extracolsep{0.03em}}l@{\extracolsep{0.03em}}l@{\extracolsep{0.03em}}l@{\extracolsep{0.03em}}}  			
			\multicolumn{2}{r}{\includegraphics[width=0.28\linewidth]{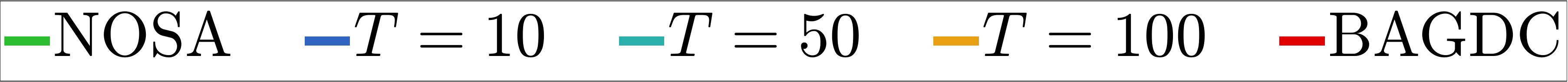}}&
			\multicolumn{2}{r}{\includegraphics[width=0.28\linewidth]{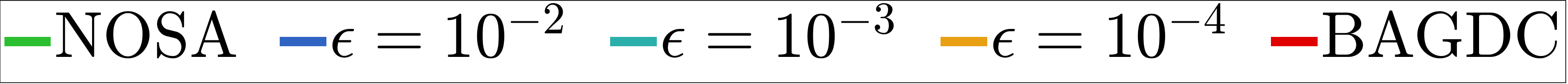}}&
			\multicolumn{2}{r}{\includegraphics[width=0.28\linewidth]{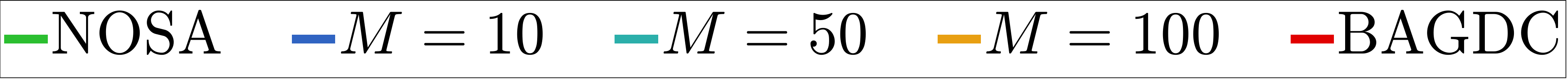}}\\
			\includegraphics[height=0.09\textheight,width=0.163\linewidth]{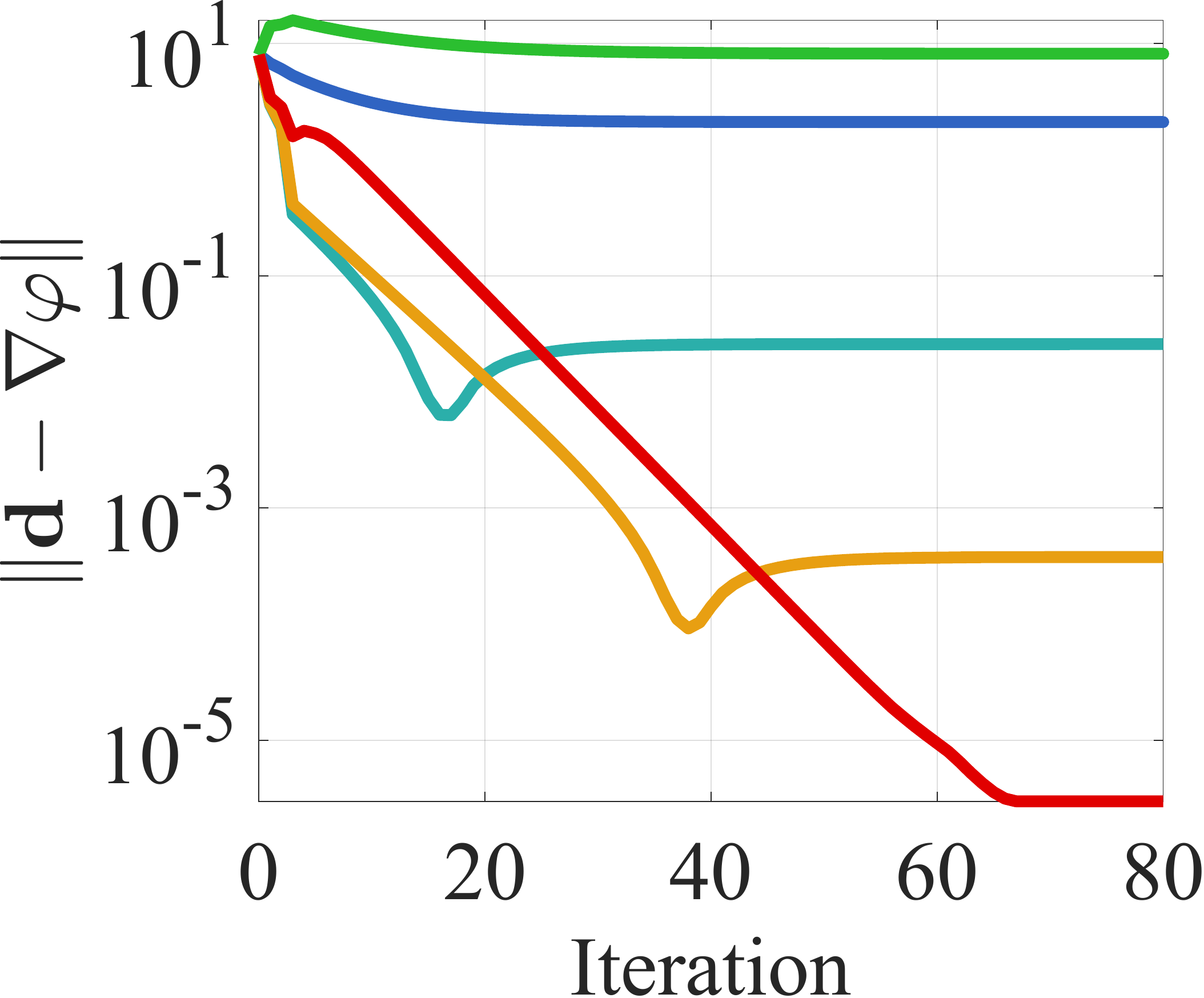}&
			\includegraphics[height=0.09\textheight,width=0.163\linewidth]{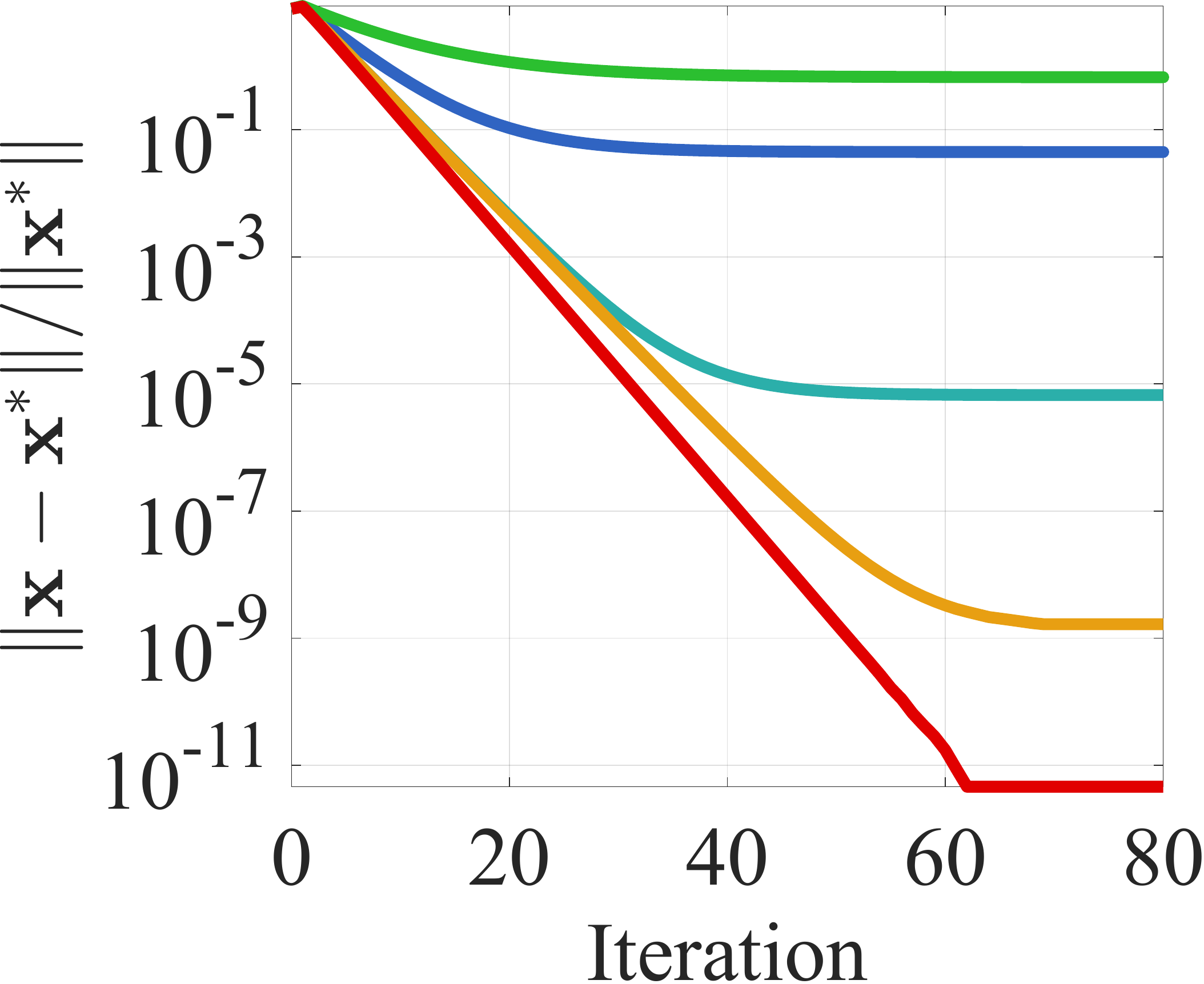}&
			\includegraphics[height=0.09\textheight,width=0.163\linewidth]{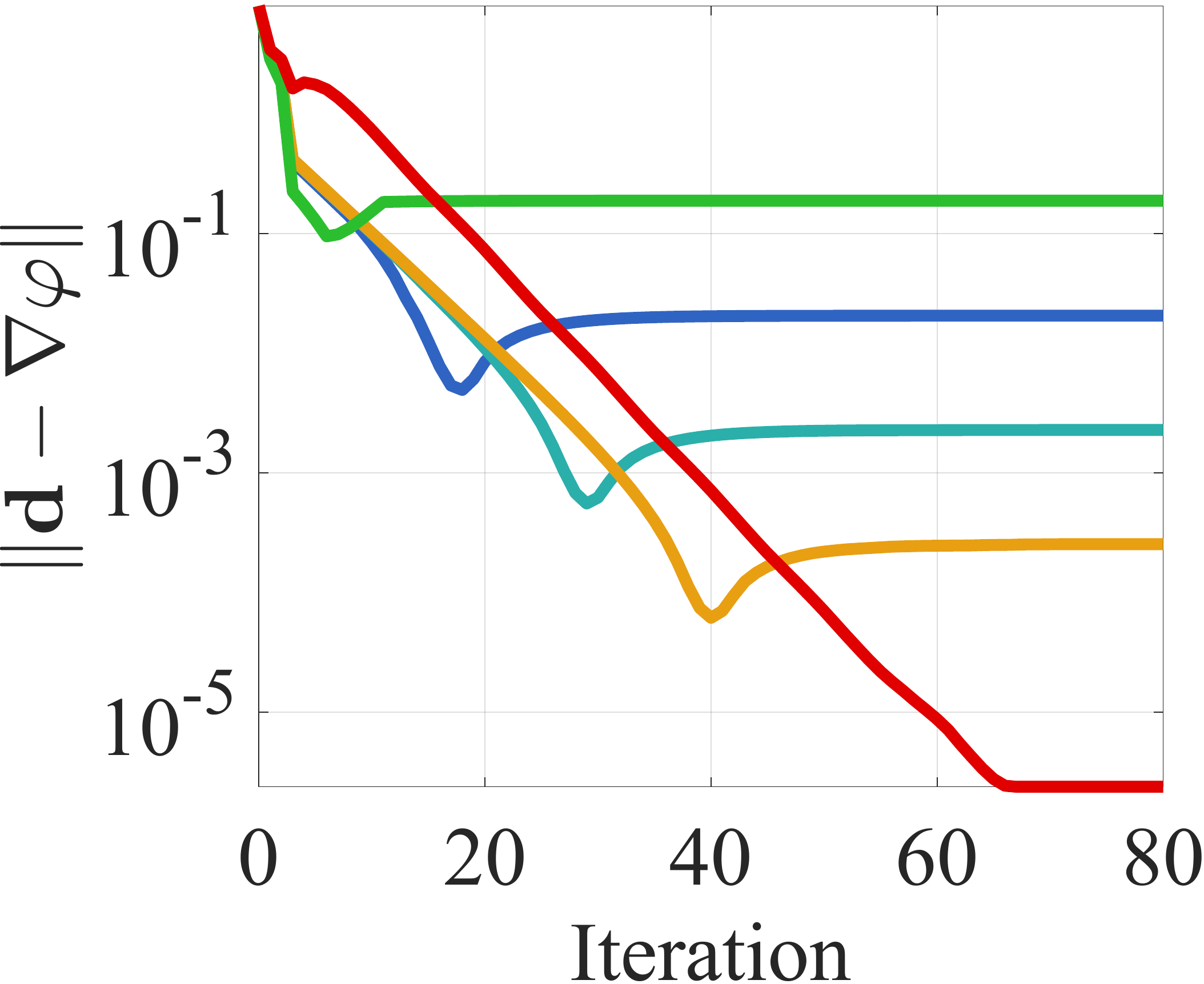}&
			\includegraphics[height=0.09\textheight,width=0.163\linewidth]{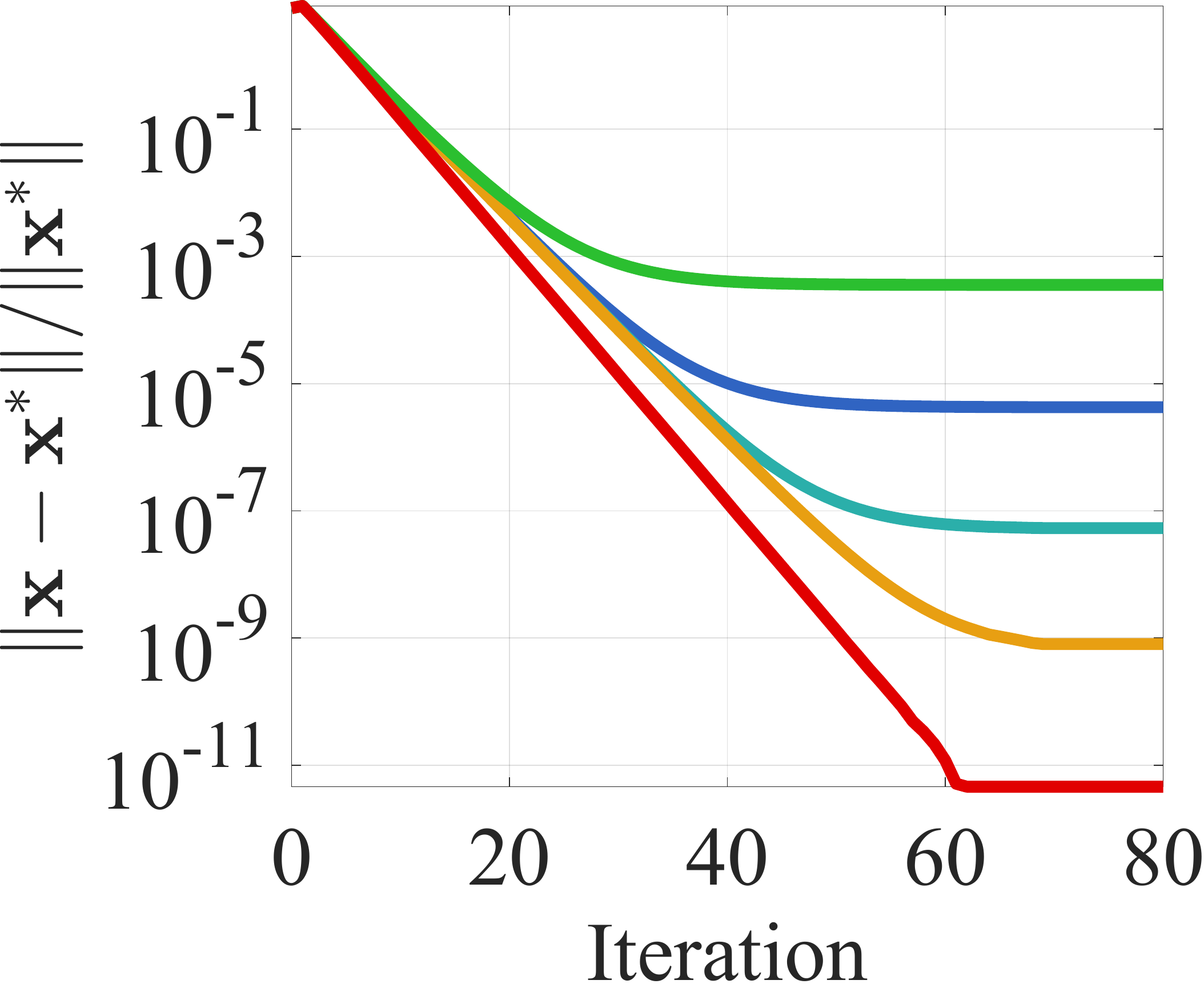}&
			\includegraphics[height=0.09\textheight,width=0.163\linewidth]{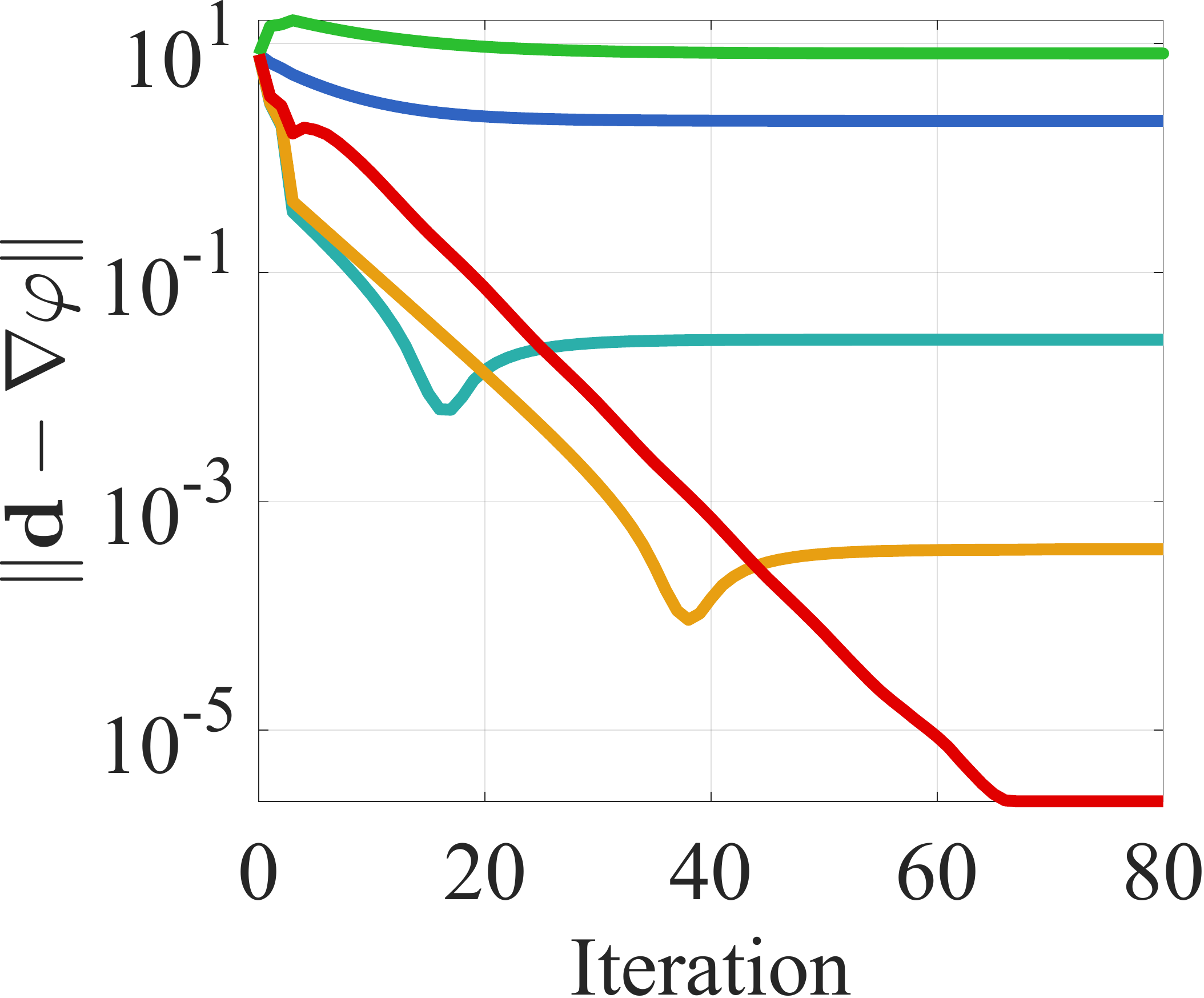}&
			\includegraphics[height=0.09\textheight,width=0.163\linewidth]{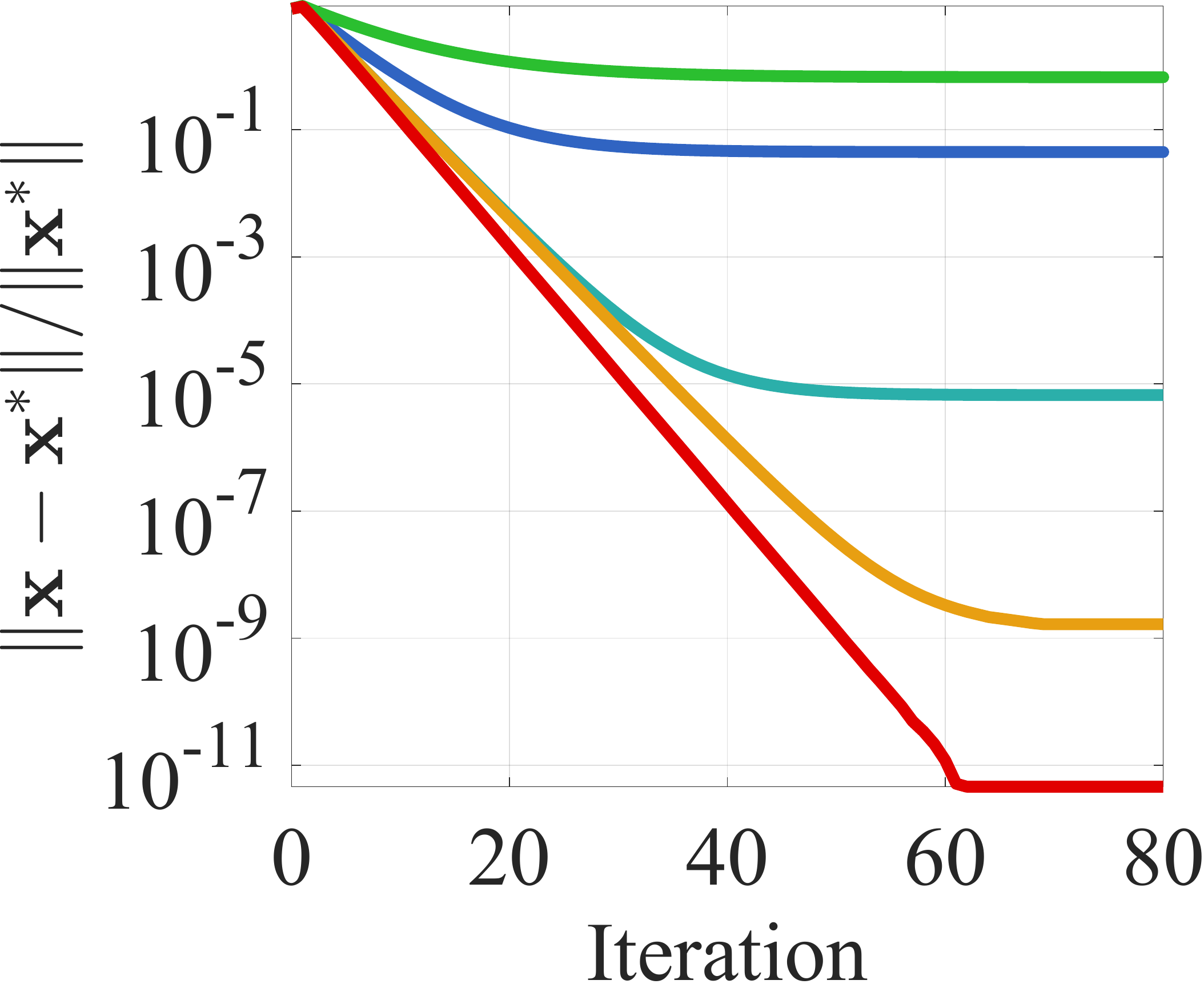}
			\\
			\multicolumn{2}{c}{\footnotesize (a) RHG}&\multicolumn{2}{c}{\footnotesize (b) CG}&\multicolumn{2}{c}{\footnotesize (c) NS}\\
		\end{tabular}
	}	 	
	
	\caption{The convergence curves of our BAGDC and various methods (RHG, CG and NS) in different accuracy of LL problem, i.e., iterations $T$ for RHG, error $\epsilon$ for CG, sequence lengths $M$ for NS. Our BAGDC method with only one step for LL problem achieves the best convergence results.}
	\label{fig:inner}
	\vspace{-0.2cm}
\end{figure}

\textbf{Solving LL with Different Accuracies vs BAGDC}.
We next show that BAGDC effectively  decouples  the dependence on the accuracy of the solution to the LL problem when applied to different BLO methods. We use the inner iterations $T$ for RHG,  error $\epsilon$ for CG and sequence lengths $M$ for NS as indicators of the accuracy of LL problem. 
In Figure \ref{fig:inner}, we compare the convergence of existing methods for different LL problem solving accuracies and existing methods with our BAGDC. We show the error of direction of descent $\Vert\mathbf{d}-\nabla F\Vert$ to illustrate that an incorrect estimate of $\v$ for the existing methods leads to the wrong direction of descent, and $\Vert\x-\x^*\Vert/\Vert\x^*\Vert$ to illustrate the effect of this on the final convergence result.
It can be seen that  accelerating existing method by directly simplifying the LL iterations will reduce the quality of the solution. By introducing $\v$ as a correction in BAGDC, our method can continuously improve accuracy for these BLOs with alternating solving UL and LL problem.

%
%Firstly, we note that all the existing methods can show the convergence trend to the best at the initial stage of the iteration, but they soon stop further convergence. The final accuracy can only be further improved by increasing the computational accuracy of the LL problem. We attribute this to the fact that the final solution accuracy of the existing methods depends entirely on the solution accuracy of the LL problem. In the earlier stages of training, the descent direction with error still can solve the problem at a coarse granularity. However, as accuracy gradually increases in the later stages of training, the final accuracy is limited by the low accuracy of the LL problem. Thanks to the correction of the multiplier variable $\v$, our BAGDC method can continuously approximate the direction of descent with each iteration of $\x$, allowing the accuracy to improve with the total number of iterations. This verifies that our method eventually obtains a reasonable gradient estimate even when the LL problem is not computed precisely at each step. 

\begin{figure}[!htbp]
	\centering
	\begin{tabular}{c@{\extracolsep{0.1em}}c@{\extracolsep{0.1em}}c@{\extracolsep{0.1em}}c@{\extracolsep{0.1em}}}
		\includegraphics[width=0.20\linewidth]{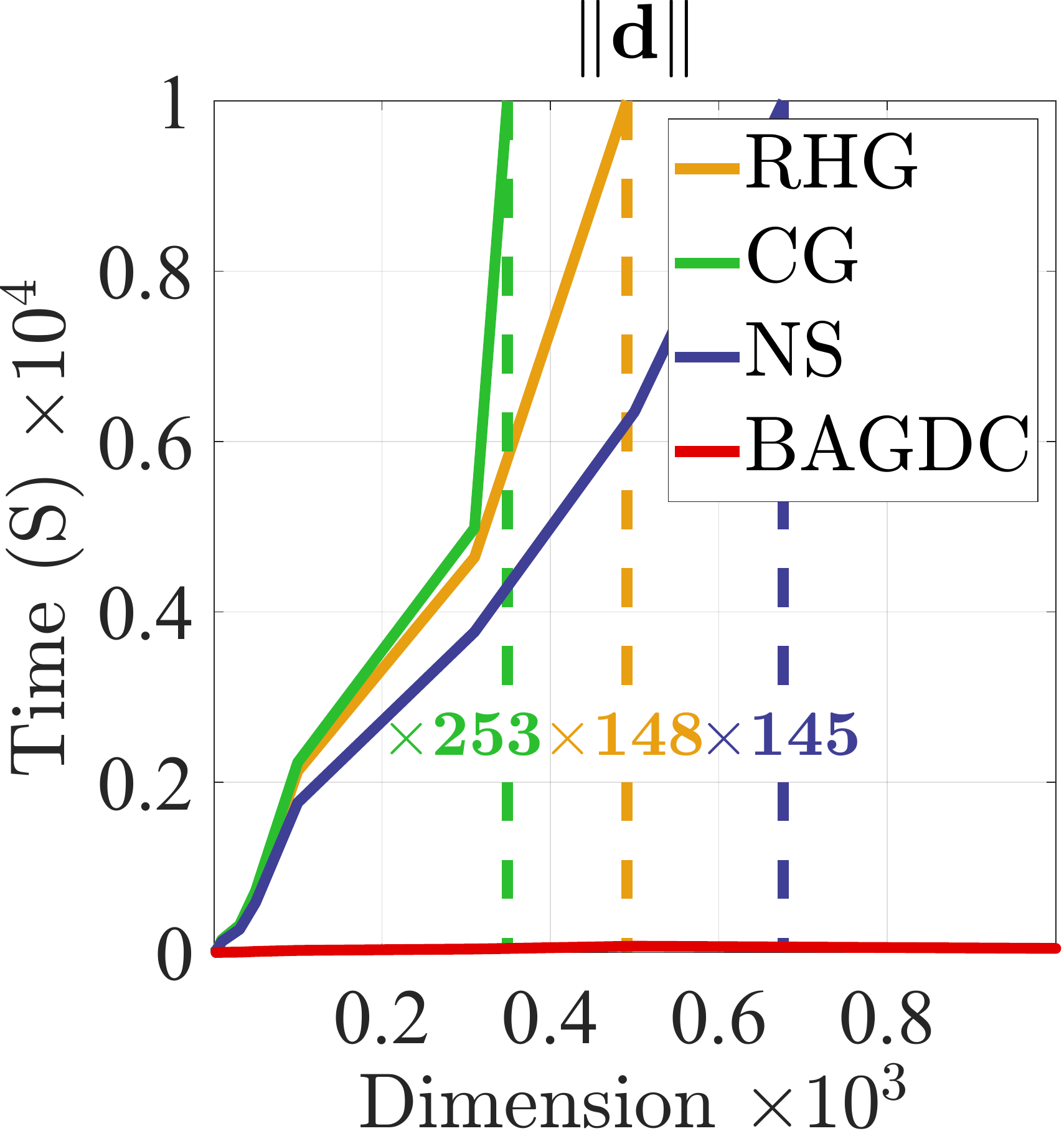}&
		\includegraphics[width=0.20\linewidth]{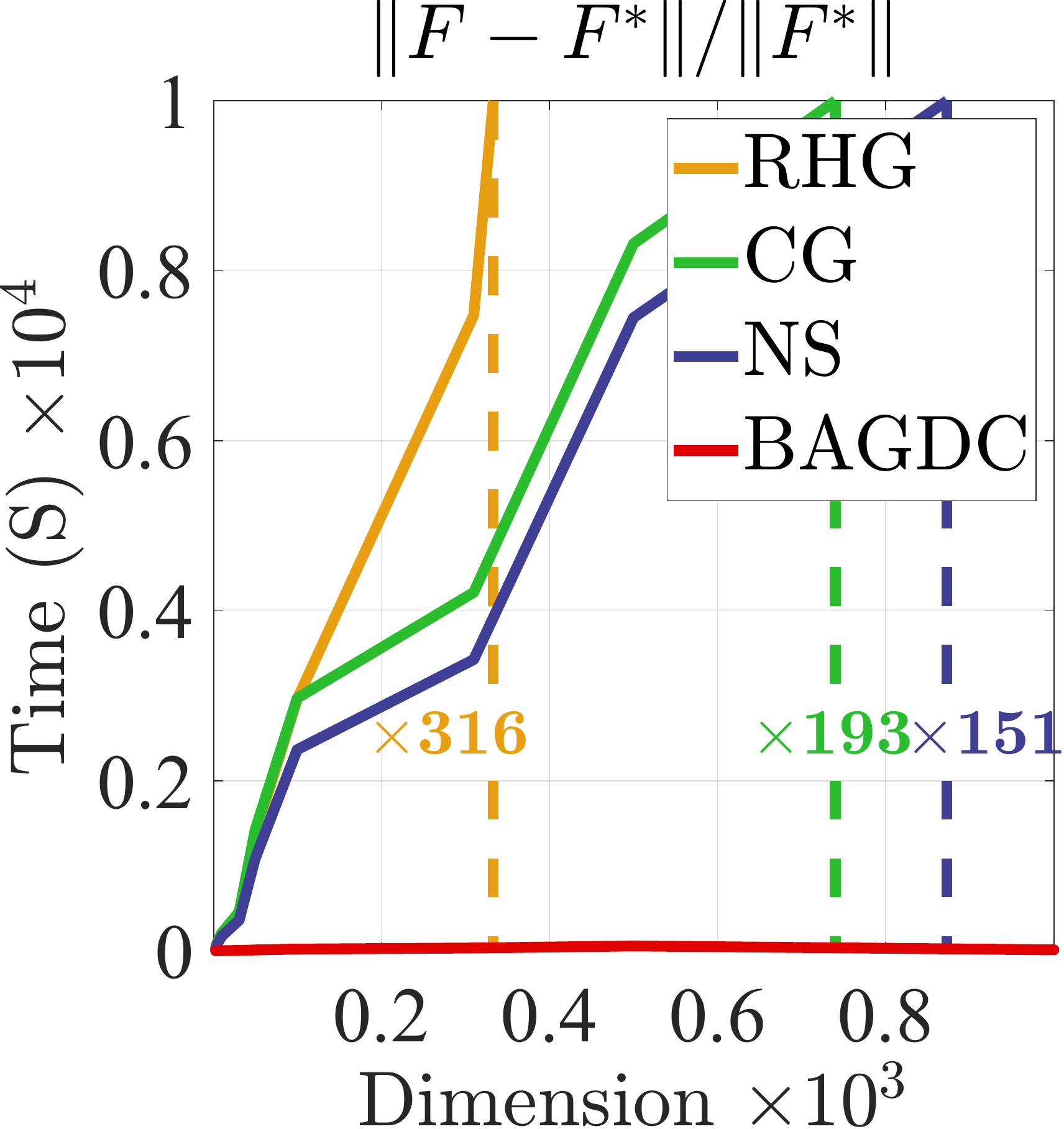}&
		\includegraphics[width=0.20\linewidth]{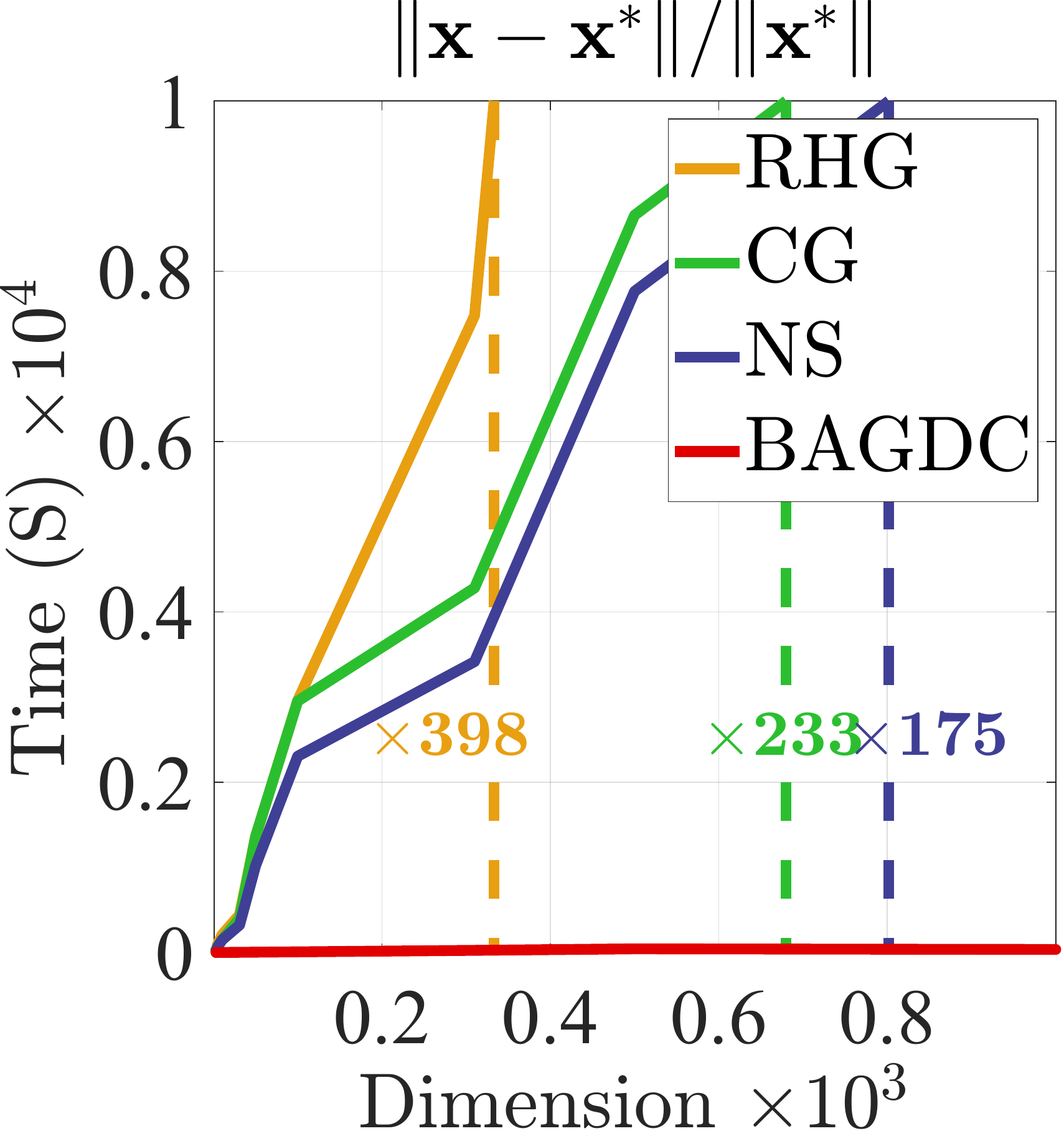}&
		\includegraphics[width=0.37\linewidth,trim=0 -34 0 0,clip]{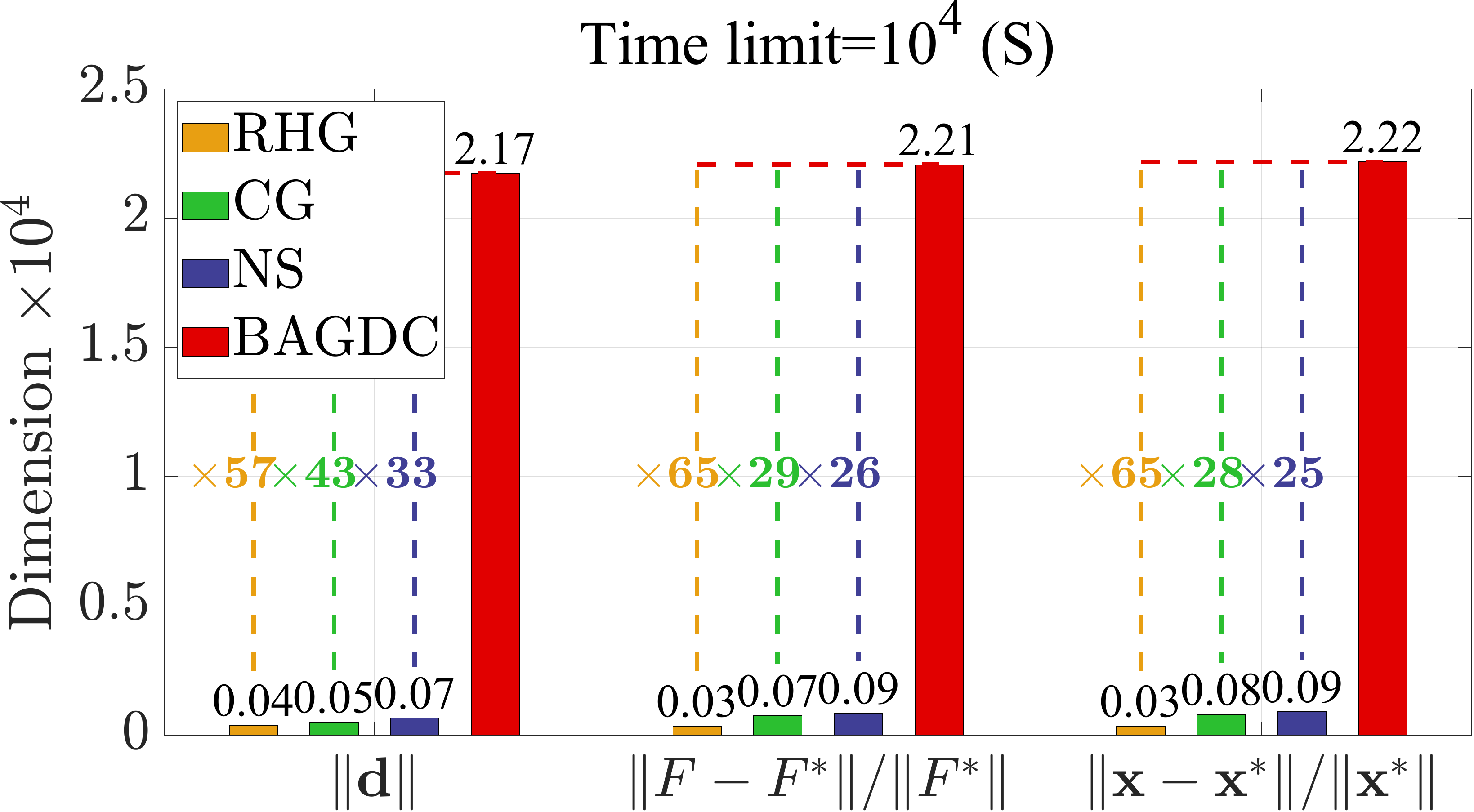}\\
		\multicolumn{3}{c}{\footnotesize (a) Convergence time}&\multicolumn{1}{c}{\footnotesize\quad (b) Convergence dimension}\\
		%	 						\multicolumn{1}{c}{$\Vert\nabla F\Vert$}&\multicolumn{1}{c}{$\Vert F-F^*\Vert/\Vert F^*\Vert$}
		%	&\multicolumn{1}{c}{$\Vert\x-\x^*\Vert/\Vert \x^*\Vert$}&\\
	\end{tabular}
	\caption{Computational efficiency comparison of existing BLOs (RHG, CG and NS) and ours BAGDC. (a) The time required for different convergence metrics in different dimensions. BAGDC reduces the time required compared to existing methods. (b) The problem dimension when reached the time limit. BAGDC greatly increases the dimension of the problem to reach the time limit.}
	\label{fig:dim}
	\vspace{-0.2cm}
\end{figure}

\textbf{Significant Speed Advantage of BAGDC.} We then show the improvement in computational efficiency of BAGDC, the extremely simpler implementation of BLO. For existing methods, the main part of the computational complexity is the large number of matrix-vector products in nested iteration. Therefore by changing the dimension of the matrix $\A$, we can easily extend Eq~(\ref{eq:counter}) to the high dimensional case to observe the computational efficiency.
Here we show the convergence in different dimensions on Eq.~(\ref{eq:counter}) with limited total time $10^4$ s. 
Figure~\ref{fig:dim} (a) shows the convergence time of different indicator $\Vert\mathbf{d}\Vert\leq 10^{-3}$, $\Vert F-F^*\Vert/\Vert F^*\Vert\leq 10^{-3}$ and $\Vert\x-\x^*\Vert/\Vert\x^*\Vert\leq 10^{-4}$. 
Since each step usually contain only one matrix vector product, the alternating optimization BAGDC saves a lot of time compared to the multi-step per iteration of nested optimization. Thus our method requires much less time compared to other methods of all scales. Figure~\ref{fig:dim}(b) illustrates the limiting problem dimension that existing methods can achieve when reaching the time limit. BAGDC does not contain repeated nested matrix-vector products allowing our method to be solved on problems of higher dimensionality.

\begin{figure}[!htbp]
	\centering
	\setlength{\tabcolsep}{1.5mm}{
		\begin{tabular}{c@{\extracolsep{0.1em}}c@{\extracolsep{0.1em}}c@{\extracolsep{0.1em}}c@{\extracolsep{0.1em}}c@{\extracolsep{0.1em}}c@{\extracolsep{0.1em}}c@{\extracolsep{0.1em}}c@{\extracolsep{0.1em}}}
			\includegraphics[height=0.11\textheight,width=0.19\linewidth]{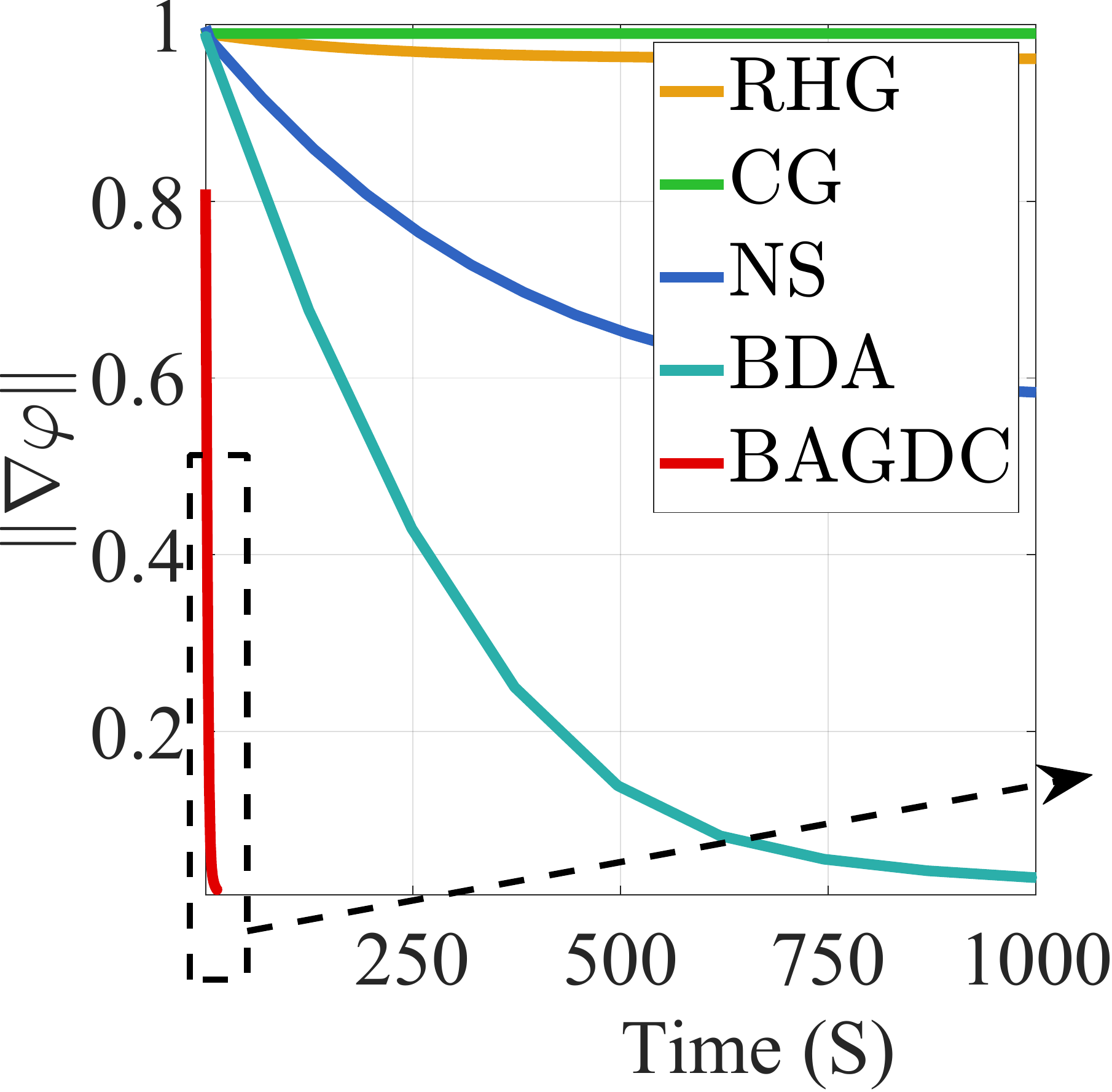}&
			\includegraphics[height=0.112\textheight,width=0.05\linewidth,trim=0 -55 0 0,clip]{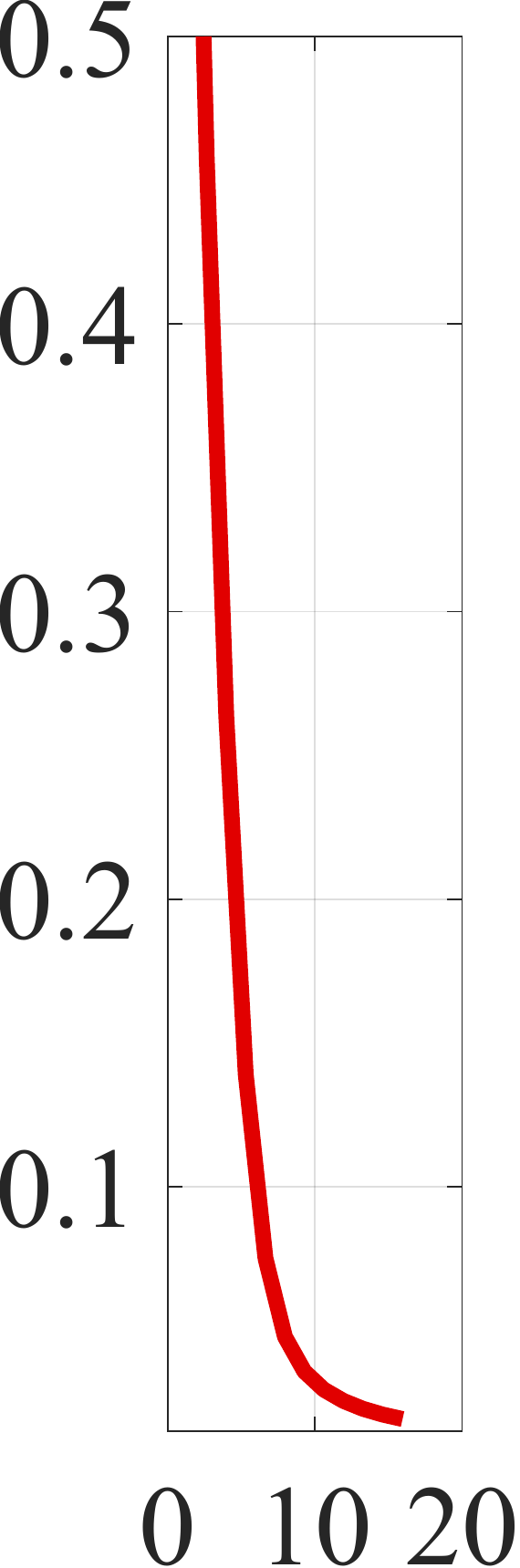}&
			\includegraphics[height=0.112\textheight,width=0.19\linewidth]{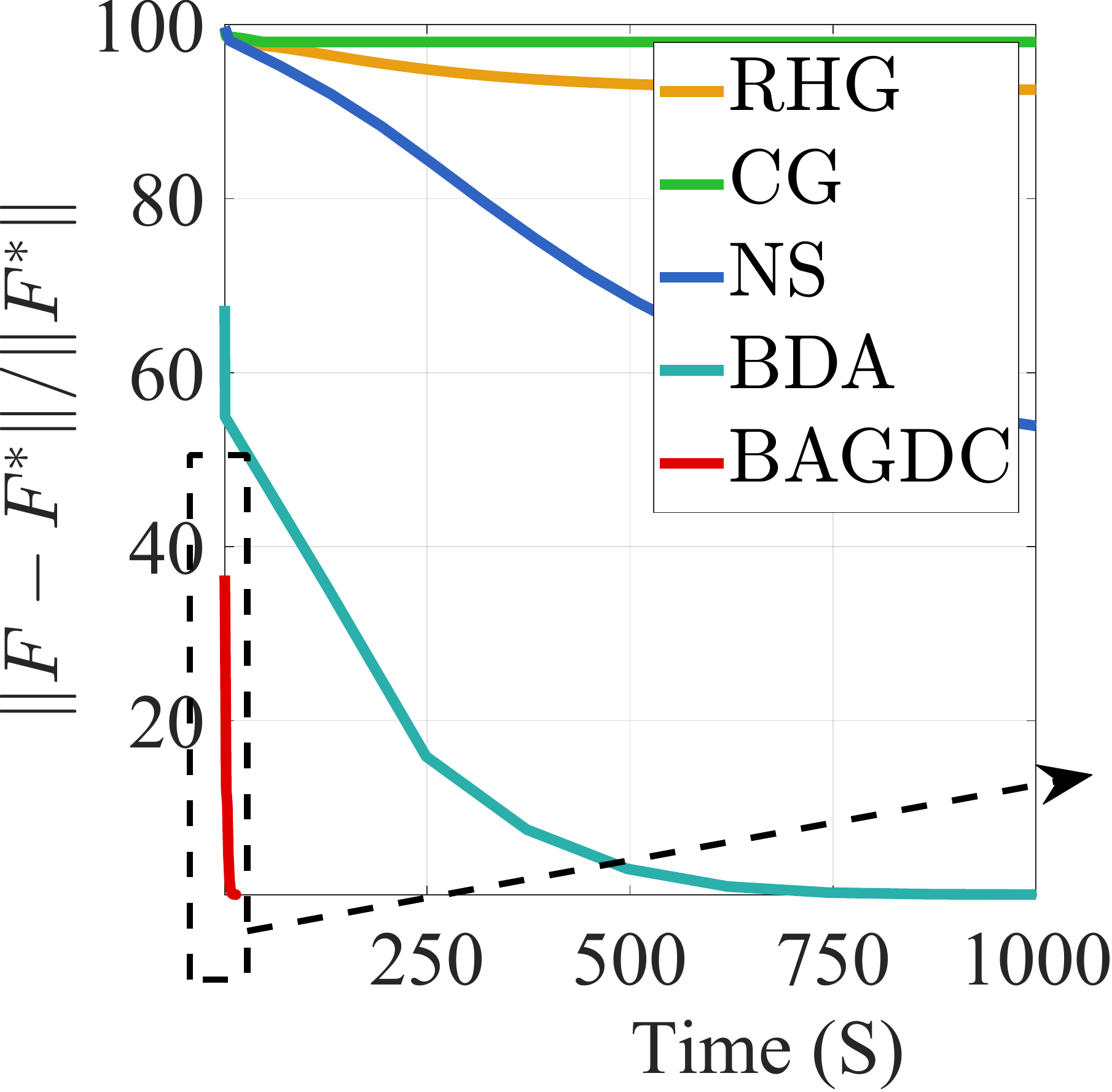}&
			\includegraphics[height=0.112\textheight,width=0.05\linewidth,trim=0 -55 0 0,clip]{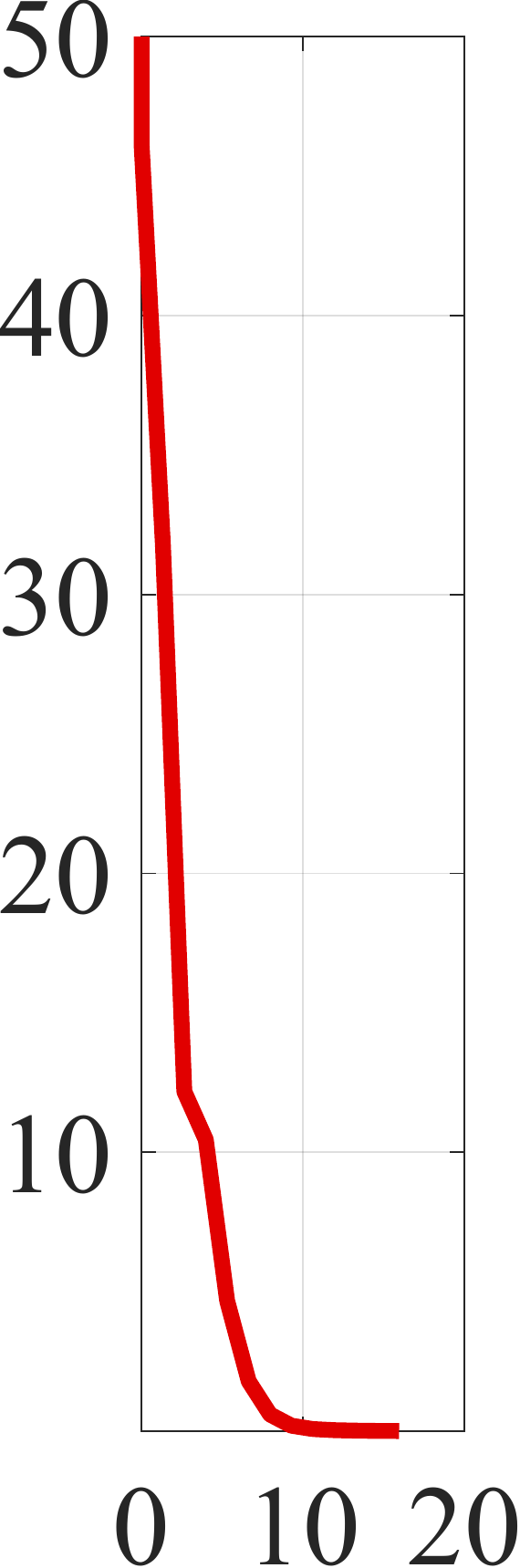}&
			\includegraphics[height=0.11\textheight,width=0.19\linewidth]{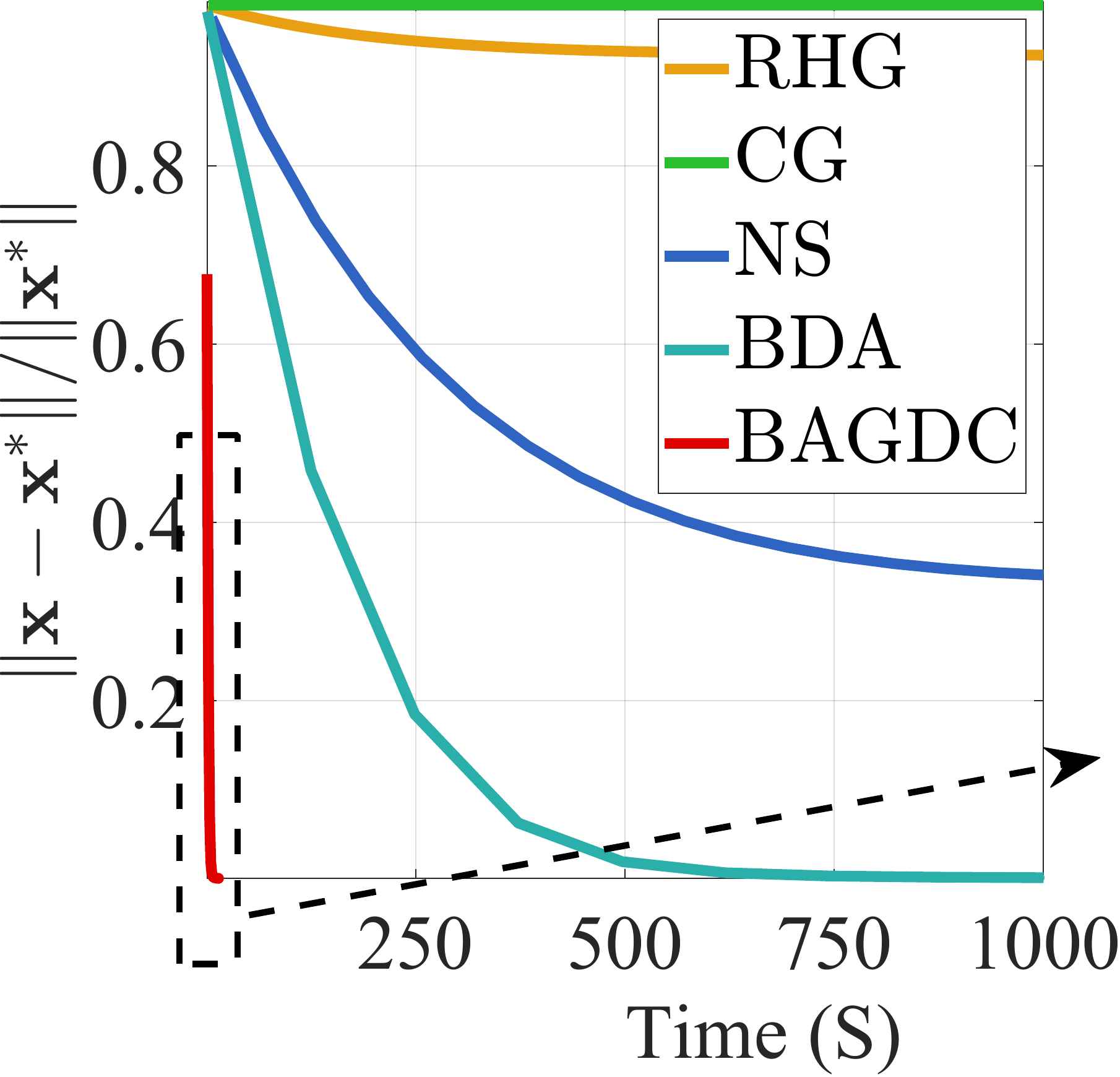}&
			\includegraphics[height=0.112\textheight,width=0.052\linewidth,trim=0 -55 0 0,clip]{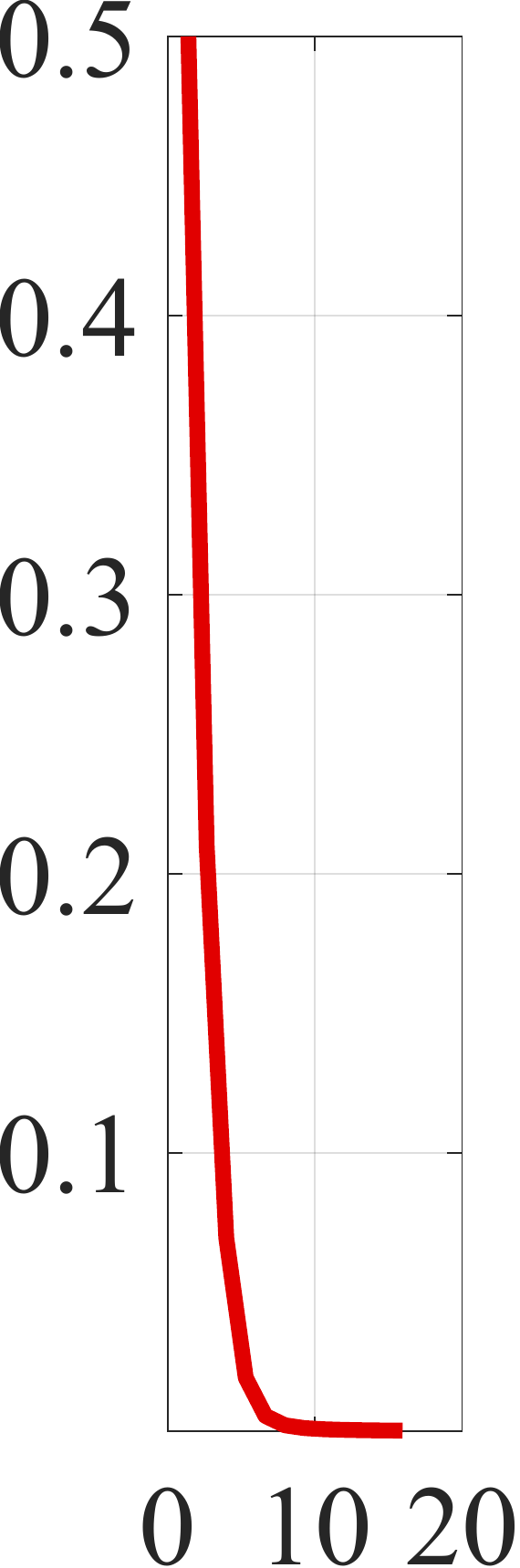}&
			\includegraphics[height=0.112\textheight,width=0.19\linewidth]{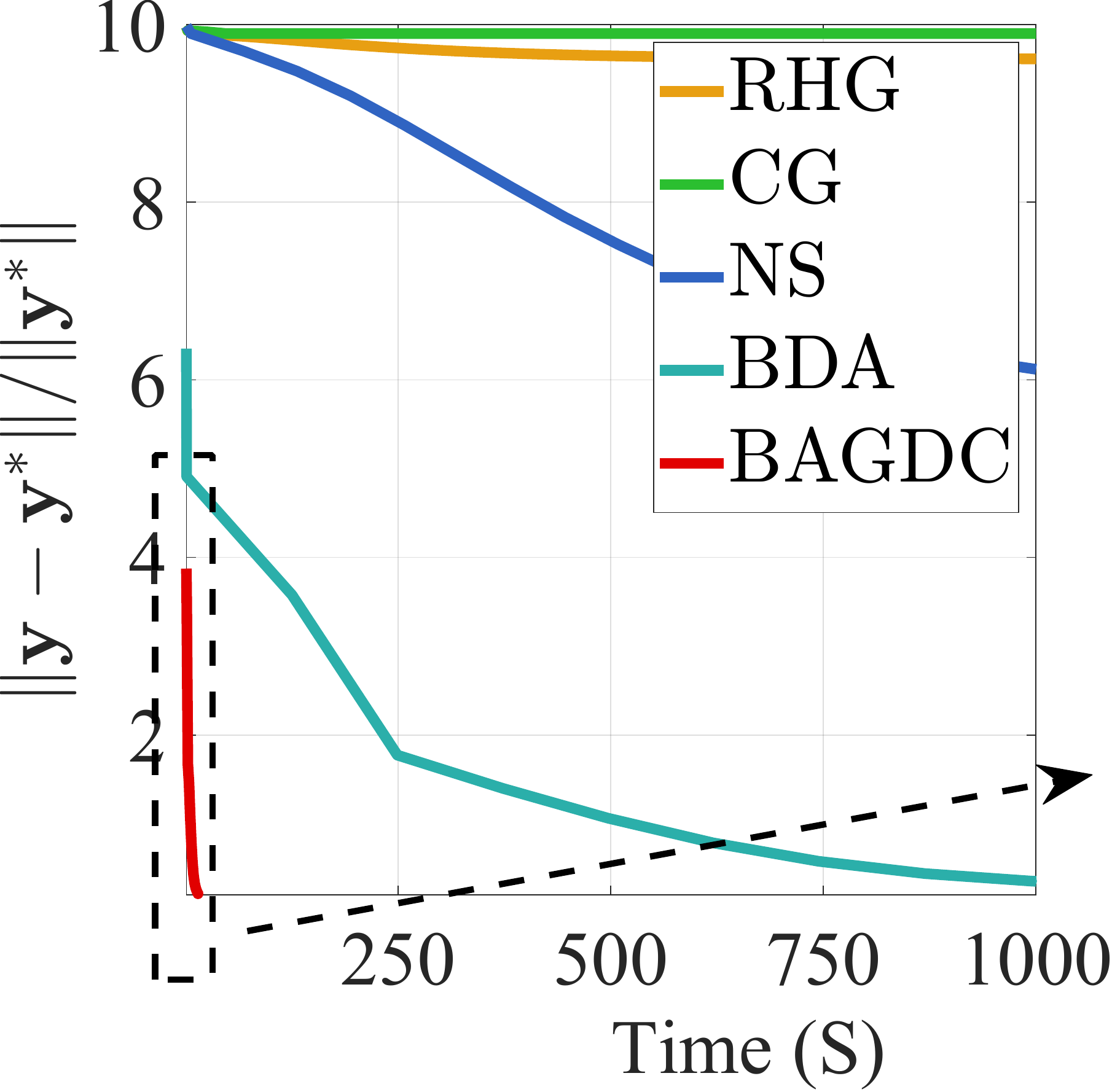}&
			\includegraphics[height=0.112\textheight,width=0.04\linewidth,trim=0 -55 0 0,clip]{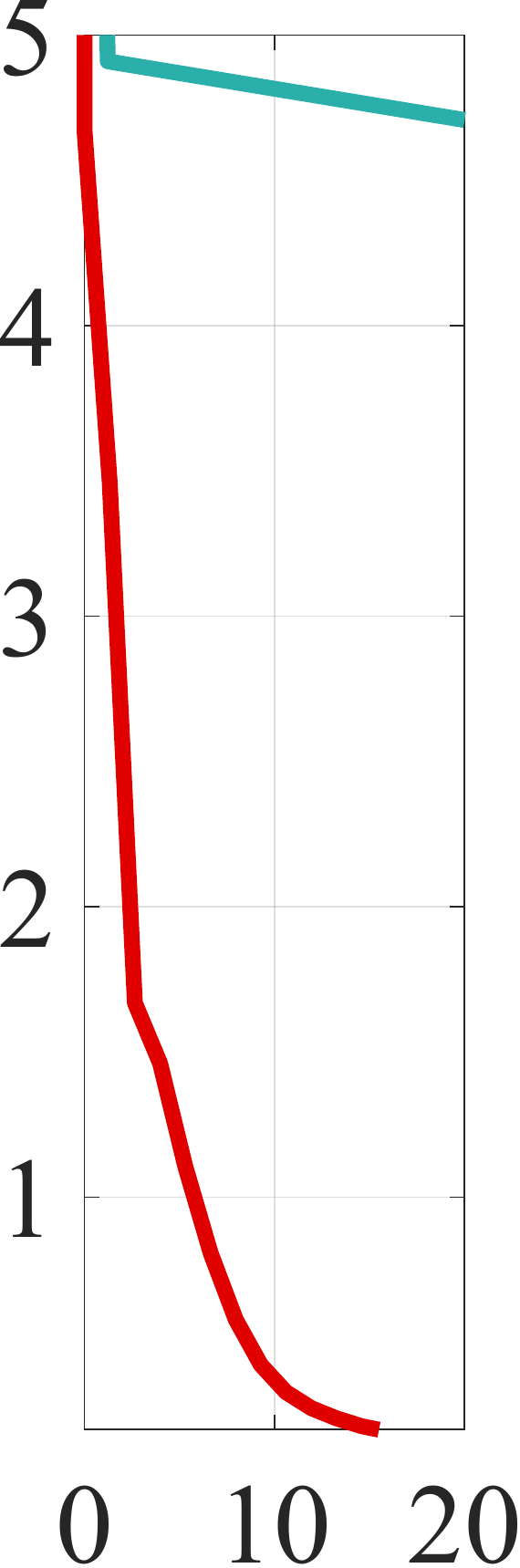}
			\\
		\end{tabular}
	}	 
	\caption{Convergence curves of existing methods under BLO beyond LLS. BAGDC has a very high convergence speed with guaranteed correct convergence because of the much simpler iteration.} 
	\label{fig:BDA}
	\vspace{-0.4cm}
\end{figure}

\textbf{Solving BLOs with multiple LL minimizers.} We next verify that BAGDC can also solve the BLOs with multiple LL minimizers  as mentioned in Theorem~\ref{thmconvex}. We use the example without LLS given in Supplemental Material.
% $F=1/2\left\| \x-\y_{2}\right\|^{2}+1/2\left\|\y_{1}-\mathbf{e}\right\|^{2}$ and $f=1/2\left\|\y_{1}\right\|^{2}-\x^{T} \y_{1}$, 
%where $\x \in \mathbb{R}^{n}$, $\y=\left(\y_{1}, \y_{2}\right) \in  \mathbb{R}^{2 n}$. The optimal solution is $\mathbf{e}\cdot\left(1, 1,  1\right)$. Note that unlike the single LL minimizers case, here $\y^* =\arg\min_{\y}F(\x,\y), \text{s.t.} \nabla_{\y}f(\x,\y)=0$.
Figure~\ref{fig:BDA} compares the convergence curves with BLO methods. As we can see, only BAGDC and BDA are able to converge to the optimum point, the other methods only converge to an incorrect stability point. Thanks to the much simpler iteration format resulting from alternating optimisation, BAGDC converge much faster than BDA.

\subsection{Real-world Applications}
In this section, we verify the performance of our BAGDC on real data and high-dimensional  problems. We show two widely used applications of BLOs, data hyper-cleaning, and few-shot learning.

\begin{wrapfigure}{r}{6.4cm}
	\vspace{-0.2cm}
	\begin{minipage}[r]{\linewidth}
		\centering 
		\makeatletter\def\@captype{figure}\makeatother 
		
		\includegraphics[width=0.48\linewidth]{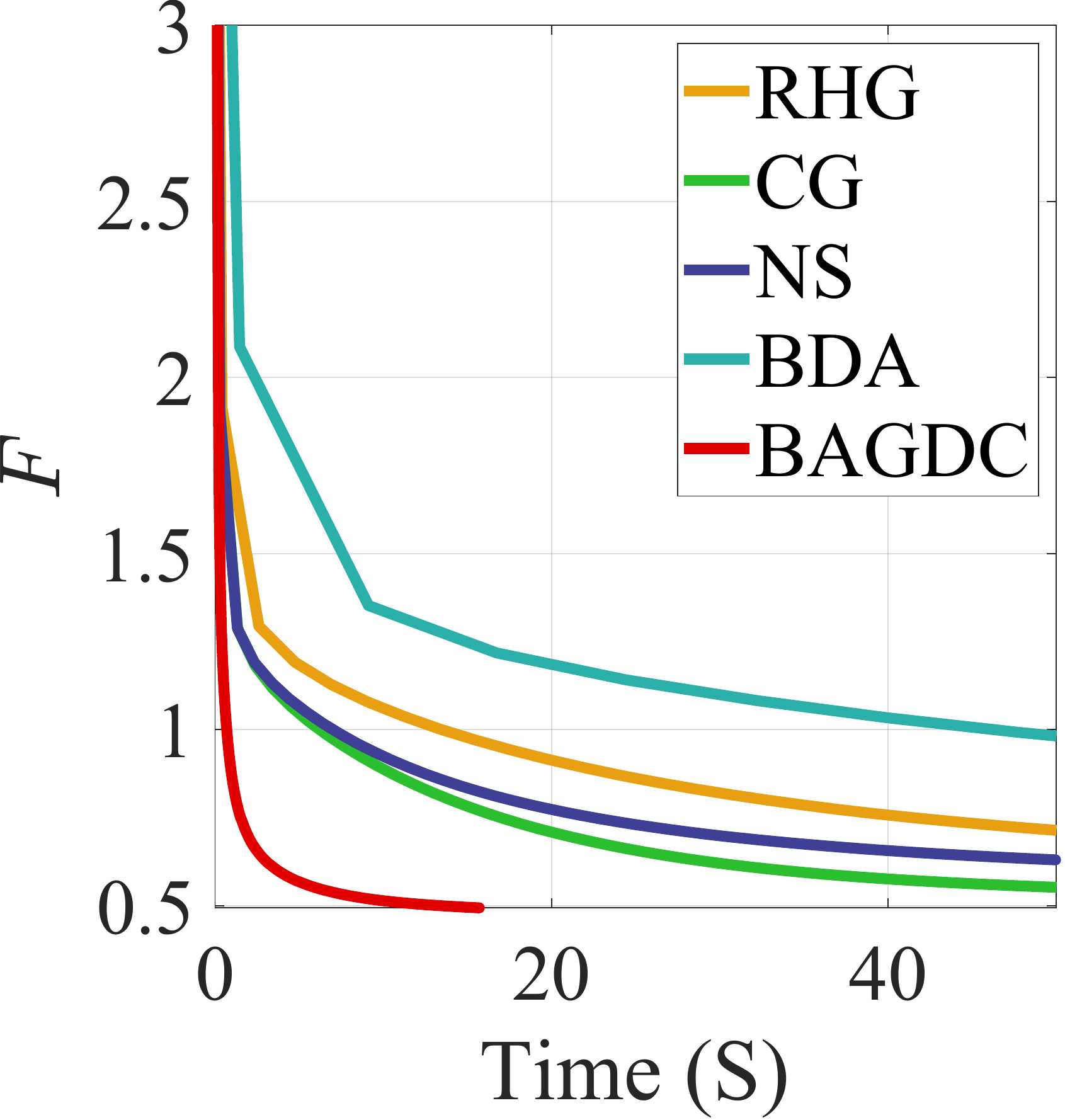}
		\includegraphics[width=0.48\linewidth]{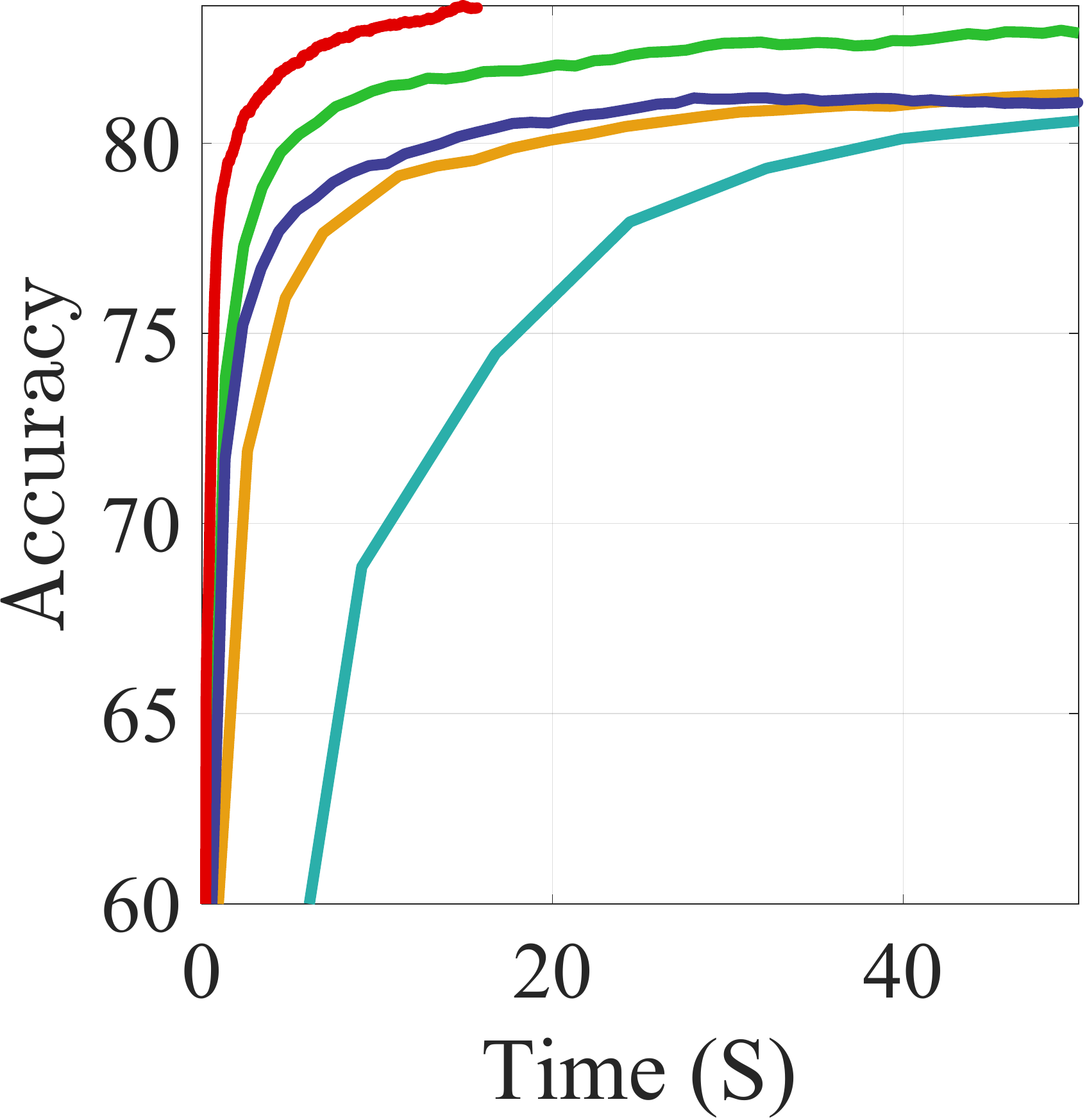}
		
		\caption{Comparison of the validation loss $F$ and accuracy for hyper-cleaning. BAGDC dramatically increases the  training speed.} 		\label{fig:valloss}%
	\end{minipage}
	\vspace{-0.4cm}
\end{wrapfigure}
\textbf{Data Hyper-cleaning.} Data hyper-cleaning is a widely used application for bi-level optimization. Assuming that some of the labels in our dataset are contaminated, data hyper-cleaning aims to reduce the impact of incorrect samples by adding hyper-parameters to them. 
In the experiment, we follow the settings in \cite{liu2020generic} and based on FashionMNIST datasets. 
To demonstrate the advantage of our method in terms of computational efficiency, we show in the left part of Table~\ref{tab:hypercleaning} the F1 scores and the time required for different methods to achieve the same accuracy. BAGDC achieves the same accuracy and similar F1 scores while reducing the time requirement to an order of magnitude.
Figure~\ref{fig:valloss} shows the validation loss for the different methods. It can be seen that our method has the fastest convergence rate and can still be the fastest at higher accuracies after quickly achieving  $81\%$ accuracy rate.

\textbf{Few-shot Learning.} We then test the application in higher dimensions to verify the computational efficiency of our method in  few-shot learning. Few-shot learning (i.e., N-way M-shot) is a multi-task N-way classification and it aims to learn the hyperparameter such that each task can be solved only with M training samples. In the experiment, we follow the settings in \cite{liu2020generic} and based on DoubleMNIST datasets~\cite{sitzmann2020metasdf}.  In the right part of Table~\ref{tab:hypercleaning} , we shows the time that existing methods achieve the same accuracy (95\% for 5-way and 85\% for-20 way). Our method achieves the same accuracy rate while significantly reducing the time required. Therefore, our alternating optimization method can be effectively applied to large-scale neural networks.

\begin{table}[h]
	\vspace{-0.3cm} 
	\centering 
	\makeatletter\def\@captype{table}\makeatother \caption{Comparison of the results of existing methods for hyper-cleaning and few-shot learning. We compared the time and F1 score for hyper-cleaning when achieving almost the same accuracy ($81\%$) and  the time for few-shot learning ($95\%$ for 5-way and $85\%$ for 20-way). BAGDC has the highest computational efficiency while training to the same accuracy rate.} 
	\setlength{\tabcolsep}{2.4mm}{
		\begin{tabular}{cccccccccc}
			\toprule
			\multirow{2}{*}{Method}&\multicolumn{3}{c}{Hyper-cleaning}&&&\multicolumn{2}{c}{ 5-way} & \multicolumn{2}{c}{20-way}\\
			&\multirow{1}{*}{Acc.}&\multirow{1}{*}{F1 score}&\multirow{1}{*}{Time}&&&Acc. &Time  &Acc.&Time \\
			\midrule
			RHG&81.07&87.10&41.47&&&  95.05  &  1614.38  & 85.06   & 2645.27\\
			BDA&81.09&88.13&64.84&&&  95.07  &  1864.90  & 85.09  & 3017.33\\
			CG&81.05&87.85&19.72&&&  95.08  &  1056.91  &  85.00  &  2060.94\\
			NS&81.03&85.53&26.01&&&  95.07  &  1005.26  & 85.05   &  1942.02\\
			BAGDC&{81.00}&87.06&\textbf{{2.92}}&&&  95.06  &  \textbf{98.16}  &  85.11   & \textbf{381.21}\\
			\bottomrule		
		\end{tabular}
	}\label{tab:hypercleaning}%
	\vspace{-0.3cm} 
\end{table}

\section{Conclusion}

This work first designed a counter-example to demonstrate that there actually exist fundamental theoretical issues in these naive single-inner-step acceleration for existing GBLOs. Then by reformulating BLO as a constrained single-level optimization problem and introducing dual multipliers from the viewpoint of KKT condition, an extremely fast gradient method, named BAGDC, is proposed for BLOs. Our theoretical results guaranteed that BAGDC can successfully obtain non-asymptotic convergence sequences towards self-governed KKT stationary points. We also extended this result for the more challenging BLO problems, i.e., with multiple LL solutions.

\bibliographystyle{unsrtnat}
\bibliography{nips2021_nonconvex_bib}

\appendix
%
%
%\section{Appendix}
%
%
%Optionally include e\xtra information (complete proofs, additional e\xperiments and plots) in the appendi\x.
%This section will often be part of the supplemental material.

\end{document}